\documentclass[onefignum,onetabnum]{siamonline220329}

\usepackage{tcolorbox}
\usepackage{scalerel}
\usepackage{caption}
\usepackage{tikz-qtree}
\usepackage{sidecap}
\usepackage{subcaption}
\usepackage{amsfonts}
\usepackage{mathtools}
\usepackage{stmaryrd}

\usetikzlibrary{arrows,decorations.markings}
\usepackage{lipsum}
\usepackage{graphicx}
\usepackage{epstopdf}
\usepackage{tikz}
\usepackage{pgfplots}
\usepackage{amssymb}
\usepackage{xcolor}
\usepackage{amsmath}

\usepackage{algorithmic}
\ifpdf
  \DeclareGraphicsExtensions{.eps,.pdf,.png,.jpg}
\else
  \DeclareGraphicsExtensions{.eps}
\fi

\usepackage{enumitem}
\setlist[enumerate]{leftmargin=.5in}
\setlist[itemize]{leftmargin=.5in}

\newsiamremark{remark}{Remark}
\newsiamremark{hypothesis}{Hypothesis}
\crefname{hypothesis}{Hypothesis}{Hypotheses}
\newsiamthm{claim}{Claim}

\headers{Sparsification of Large Strictly Ultrametric Matrices}{E.D. Gorman and M.E. Lladser}

\title{Sparsification of Large Ultrametric Matrices: Insights into the Microbial Tree of Life
\thanks{Submitted to the editors DATE.
\funding{This work was partially funded by the NSF grant No. 1836914.}}}

\author{Evan D. Gorman\footnotemark[2]
\and Manuel E. Lladser\thanks{Department of Applied Mathematics, University of Colorado, Boulder, CO 80309, The United States
  (Corresponding author e-mail: \email{manuel.lladser@colorado.edu)}}
}

\usepackage{amsopn}

\ifpdf
\hypersetup{
  pdftitle={Sparsification of Large Ultrametric Matrices: Insights into the Microbial Tree of Life},
  pdfauthor={E. D. Gorman and M. E. Lladser}
}
\fi

\newcommand{\avg}{\text{avg}}
\newcommand{\var}{\text{var}}
\newcommand{\diag}{\text{diag}}
\newcommand{\EE}{\mathbb{E}}

\newcommand{\EPL}{\text{EPL}}
\newcommand{\IPL}{\text{IPL}}
\newcommand{\TPL}{\text{TPL}}

\newcommand{\indicator}[1][\ ]{[\![ {#1} ]\!]} 
\newcommand{\PP}{\mathbb{P}}
\newcommand{\RR}{\mathbb{R}}
\newcommand{\TT}{\mathbb{T}}
\newcommand{\VV}{\mathbb{V}}
\newcommand{\Onebf}{\mathbf{1}}

\begin{document}

\maketitle

\begin{abstract}
Ultrametric matrices have a rich structure that is not apparent from their definition.  Notably, the subclass of strictly ultrametric matrices are covariance matrices of certain weighted rooted binary trees.  In applications, these matrices can be large and dense, making them difficult to store and handle.  In this manuscript, we exploit the underlying tree structure of these matrices to sparsify them via a similarity transformation based on Haar-like wavelets.  We show that, with overwhelmingly high probability, only an asymptotically negligible fraction of the off-diagonal entries in random but large strictly ultrametric matrices remain non-zero after the transformation; and develop a fast algorithm to compress such matrices directly from their tree representation.  We also identify the subclass of matrices diagonalized by the wavelets and supply a sufficient condition to approximate the spectrum of strictly ultrametric matrices outside this subclass. 
Our methods give computational access to a covariance model of the microbiologists' Tree of Life, which was previously inaccessible due to its size, and motivate defining a new but wavelet-based phylogenetic $\beta$-diversity metric.  Applying this metric to a metagenomic dataset demonstrates that it can provide novel insight into noisy high-dimensional samples and localize speciation events that may be most important in determining relationships between environmental factors and microbial composition.
\end{abstract}



\section{Introduction}

Ultrametric matrices appear across many domains of mathematics and science. They comprise an important class of matrices called inverse-M matrices \cite{DelMarSan14} and are a key object of study in potential theory and Markov Chains \cite{delmar96}. In scientific applications, ultrametric matrices act as covariance models in phylogenetic comparative analysis \cite{pur11}, network inference \cite{jita19} and energy models in statistical physics \cite{cap87}. Further hinting at the pervasiveness of ultrametric matrices in modern data science, recent work has shown that the matrix of normalized Euclidean distances between points in some random subsets of $\RR^d$ converge in probability to an ultrametric matrix as the dimension $d$ tends to infinity~\cite{Zub14,Zub17}. 

In many applications, the underlying ultrametric matrix can be massive, potentially too large to store in computer memory. This raises many challenges in the analysis and application of such matrices. However, if a sparse representation of a matrix can be found, many otherwise impossible tasks become computationally feasible. Examples of the latter include solving linear equations, matrix factorizations,  eigenvalue decompositions, and principal component analysis (PCA). 

In this paper we focus on the subclass of strictly ultrametric matrices, which is provided in the following definition.

In what remains of this manuscript, $n\ge1$ denotes an integer. Define $[n]:=\{1,\ldots,n\}$. In addition, vectors and sometimes functions are represented as column ones. The transpose of a vector or matrix $A$ is denoted $A'$.

\begin{tcolorbox}
\begin{definition}[\cite{varnab}]
A matrix $S\in\RR^{n\times n}$ is ultrametric if it is symmetric with nonnegative entries and $S(i,j)\geq\min\{S(i,k),S(k,j)\}$ for all $i,j,k\in[n]$. If in addition $S(i,i)>\max\{S(i,t):t\ne i\}$, for all $i\in[n]$, $S$ is called strictly ultrametric. For $n=1$, the last inequality is replaced with $S(i,i)>0$.
\end{definition}
\end{tcolorbox}

Ultrametric matrices have rich properties that are not made evident by their definition~\cite{DelMarSan14}. In particular, if $S$ is strictly ultrametric then it is positive definite (hence invertible), $S^{-1}$ is strictly diagonally dominant with non-positive off-diagonal entries, and $S(i,j)=0$ if and only if $S^{-1}(i,j)=0$. These properties were initially proved using probabilistic methods~\cite{MarMicSan94}. An alternative proof is based on an equivalence between strictly ultrametric matrices and a subclass of binary trees~\cite{NabVar94}. The key ingredient for this equivalence is that for $n>1$, if $S\in\RR^{n\times n}$ is symmetric with non-negative entries then it is strictly ultrametric if and only if there exists a permutation matrix $P$ and strictly ultrametric matrices $A$ and $B$ such that
\begin{equation}
P\big(S-\min(S)\,\Onebf\Onebf'\big)P'=\begin{bmatrix}
A & 0\\
0 & B\\
\end{bmatrix},
\label{ide:NabVar94}
\end{equation}
where $\min(S)$ is the smallest entry in $S$, and $\Onebf\in\RR^n$ is the column vector of ones~\cite[Proposition 2.1]{NabVar94}. Since $A$ and $B$ are of the same kind as $S$, this process may be applied recursively and the matrix $S$ encoded as a weighted rooted binary tree with special characteristics. Here we adopt a slightly different encoding to the one in~\cite{NabVar94}, which is more suitable for our purposes. The reader unfamiliar with the jargon and notation of trees may skip ahead to Section~\ref{subsec:terminology} and come back to make better sense of the construction below.

Let $S$ be a strictly ultrametric matrix of dimensions $n\times n$. We can represent $S$ as a rooted binary tree with $2n$ vertices (of which half are leaves) and hence $(2n-1)$ edges, satisfying the following definition.

\begin{tcolorbox}
\begin{definition}
An out-rooted bifurcating tree (ORB-tree) with $n$ leaves is a weighted rooted tree with the following properties: each vertex has degree $1$ or $3$; its leaf set is $[n]$ and excludes the root, which has degree 1; each edge is labeled by the subset of leaves that descend from it; and the length $\ell(e)$ of each edge $e$ is non-negative but $\ell(e)>0$ when $e$ connects a leave with its parent.
\end{definition}
\end{tcolorbox}

\begin{figure}
\centering
\begin{minipage}{0.49\textwidth}
\centering
\begin{tikzpicture}[scale=.65, node distance={15mm}, thick, main/.style = {draw, circle}] 
\node[shape=circle,draw=black,minimum size=.65cm] (A) at (0,0) {};
\node[shape=circle,draw=black,minimum size=.65cm] (B) at (0,-2) {};
\node[shape=circle,draw=black] (C) at (1.5,-4) {3};
\node[shape=circle,draw=black,minimum size=.65cm] (D) at (-1.5,-4) {};
\node[shape=circle,draw=black] (E) at (-3,-6) {1};
\node[shape=circle,draw=black] (F) at (0,-6) {2};
    \path [-] (A) edge node[left] {$\{1,2,3\}$} (B);
    \path [-] (B) edge node[left] {$\{1,2\}$} (D);
    \path [-] (B) edge node[right] {$\{3\}$} (C);
    \path [-] (D) edge node[left] {$\{1\}$} (E);
    \path [-] (D) edge node[right] {$\{2\}$} (F);
\end{tikzpicture} 
\end{minipage}
\begin{minipage}{0.49\textwidth}
\centering
\begin{tabular}{|ccc|}
\hline
$e$ & $\ell(e)$ & $\delta_e$\\
\hline
$\{1,2,3\}$ & 2 & $(1,1,1)'$\\
$\{1,2\}$ & 0 & $(1,1,0)'$\\
$\{1\}$ & 3 & $(1,0,0)'$\\
$\{2\}$ & 1 & $(0,1,0)'$\\
$\{3\}$ & 2 & $(0,0,1)'$\\
\hline
\end{tabular}
\end{minipage}
\caption{\textbf{ORB-tree and strictly ultrametric matrix correspondence.} The matrix encoding of the tree is \begin{small}$\begin{pmatrix}
5 & 2 & 2 \\
2 & 3 & 2 \\
2 & 2 & 4
\end{pmatrix}
=2\,\delta_{\{1,2,3\}}\delta_{\{1,2,3\}}'+0\,\delta_{\{1,2\}}\delta_{\{1,2\}}'+3\,\delta_{\{1\}}\delta_{\{1\}}'+\delta_{\{2\}}\delta_{\{2\}}'+2\,\delta_{\{3\}}\delta_{\{3\}}'$.\end{small}}
\label{fig:A+Tree}
\vspace{-18pt}
\end{figure}
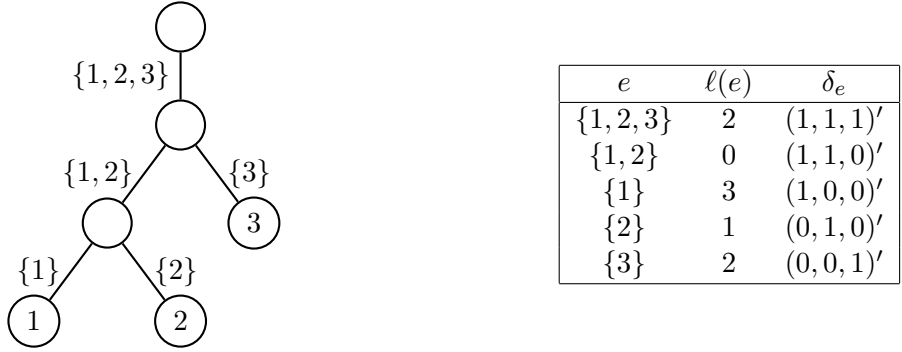

The representation of a strictly ultrametric matrix $S$ as an ORB-tree may be obtained as follows. The only edge emanating from the root is labeled as $[n]$ and defined to have length $\min(S)$. The only child of the root has two children. One child descends from an edge labeled by the rows (or columns) of $S$ associated with the matrix $A$ before applying the permutation matrix $P$ in (\ref{ide:NabVar94}). This edge has length $\min(A)$. Likewise, the other child descends from an edge labeled by the rows associated with the matrix $B$ and has length $\min(B)$. Since $A$ and $B$ are strictly ultrametric, just of smaller dimensions, the tree may be grown recursively from any descendent of the root that is not associated with a strictly ultrametric matrix of dimensions $1\times 1$. The latter represent edges that parent a leave in the ORB-tree. These edges must have a strictly positive length because $1\times 1$ strictly ultrametric matrices are strictly positive real numbers. To fix ideas see Figure~\ref{fig:A+Tree}.

For each edge $e$ in the ORB-tree define $\delta_e$ as the vector of dimension $n$ with entries $\delta_e(i)=1$ for $i\in e$ and $\delta_e(i)=0$ for $i\in[n]\setminus e$. It follows from~\cite[Theorem 2.2]{NabVar94} that
\begin{equation}
S=\sum_{e\in E}\ell(e)\,\delta_e\,\delta_e'.
\label{ide:ASumEs}
\end{equation}
Conversely, starting with any ORB-tree with edge set $E$ and length function $\ell:E\to[0,\infty)$, the matrix $\sum_{e\in E}\ell(e)\,\delta_e\,\delta_e'$ is strictly ultrametric. This representation of strictly ultrametric matrices as ORB-trees is therefore one-to-one. In fact, for $i,j\in[n]$, the entries of the matrix $S$ associated with an ORB-tree can be computed by direct inspection of tree using that
\begin{equation}
S(i,j)=\sum_{e\in[i\wedge j,\circ]}\ell(e),
\label{ide:Sij}
\end{equation}
where the $\circ$ denotes the root of the tree. We may say therefore that the entries of a strictly ultrametric matrix are indexed by the leaves of its associated ORB-tree.

We call a matrix with entries such as (\ref{ide:Sij}) the covariance matrix of the ORB-tree. This terminology is borrowed from the ecology literature where matrices like this are commonly referred to as a tree-structured or phylogenetic covariance matrices. In this setting, the leaves represent organisms, and the matrix entries denote a trait's covariance between pairs of organisms~\cite{cavsfo67}. (The term of cophenetic matrix or cophenetic distance has also been used occasionally in the hierarchical clustering literature~\cite{sar13}.)

A key tool in understanding and analyzing non-stationary and noisy continuous signals are wavelets: localized, wave-like functions.  Traditional wavelets are defined only in Euclidean spaces and have been remarkably successful in identifying multiscale structures in signals and producing sparse representations of the same. \cite{mal09} 

The Haar wavelet is among the oldest and involves averaging a signal locally at different time or space scales~\cite{graps95}. Recently, the authors of~\cite{gavnad10} extended it past continuous signals introducing the Haar-like wavelet. This new wavelet is designed for the multiscale analysis of discrete datasets equipped with a partition tree---a hierarchical structure that clusters the data into smaller subsets recursively. Due to the organization of such datasets into different tree levels (i.e., scales) and clusters (i.e., localizations), Haar-like wavelets may identify meaningful structures in data that may be impossible to distinguish otherwise---especially in noisy high dimensional datasets. 

Strictly ultrametric matrices can be fully dense; i.e all of their entries be non-zero. Nevertheless, due to the identity in equation (\ref{ide:Sij}), their entries contain much redundancy, suggesting they may be amenable to some form of compression.

This paper exploits the equivalence between strictly ultrametric matrices and ORB-trees to sparsify and hence compress the former via a change of basis. This basis is composed of the so-called Haar-like wavelets of the associated ORB-trees. The sparsification achieved by these wavelets can be substantial in large, strictly ultrametric matrices, giving computational access to matrices previously inaccessible due to their size. This can be of great value in extensive phylogenetic studies due to the interpretation of these matrices as covariance matrices of phylogenetic trees. It may also find practical applications in the context of double principal coordinate analysis, a metric of phylogenetic diversity among microbial environments.

\subsection{Paper organization}

In Section~\ref{sec:Haarlike} we introduce the Haar-like basis from \cite{gavnad10} and give a geometric interpretation of its action on ORB-trees. Then, Section~\ref{sec:Sparisfication} presents the conditions under which the Haar-like basis can be used to sparsify large, strictly ultrametric matrices. We show that the basis can substantially sparsify most large random ORB-tree's covariance matrices. We also present an algorithm for directly computing the sparsified matrix from the ORB-tree of the original matrix; in particular, without having to pre-compute the strictly ultrametric matrix from the tree. Following in Section~\ref{sec:Spectrum}, we detail the case in which the Haar-like basis diagonalizes (i.e., fully sparsifies) a strictly ultrametric matrix and provide examples of well-known tree topologies. And, in Section~\ref{sec:Approx}, we show that the conditions necessary for diagonalization can be relaxed and the Haar-like basis used to estimate eigenvalues of ORB-tree's covariance matrices.

Finally, in Section~\ref{sec:treeoflife}, we apply our methods to a covariance model of the 97\% Greengenes tree, a standard representative phylogeny of microbiologists’ Tree of Life. The sparsification opens the door for otherwise impossible tasks related to this model, such as computing the spectrum or inverse of its covariance matrix---standard tasks in phylogenetic comparative methods. We also introduce a new wavelet-based phylogenetic ($\beta$-diversity) metric corresponding to a multiscale analysis of organism abundances in microbial environments. This novel metric gives remarkably similar results to other well-known metrics on a previously studied dataset; however, it can also determine the speciation events responsible for the observed microbial compositions and quantify their respective importance.

\subsection{Paper notation and terminology}
\label{subsec:terminology} 

For real-vectors $x=(x_i)_{1\le i\le k}$ and $y=(y_i)_{1\le i\le k}$ of dimension $k$, we define $\langle x,y\rangle:=x'y=\sum_{i=1}^kx_iy_i$ and $\|x\|_2:=\sqrt{\langle x,x\rangle}$. Also, $\indicator[\cdot]$ denotes the indicator function of the proposition within the parentheses.

In our context, trees are finite undirected connected graphs without cycles.

In what remains of this manuscript, $T$ denotes an ORB-tree with $n$ leaves and branch length function $\ell:E\to\RR$. We denote the vertex and edge set of $T$ as $V$ and $E$, respectively. The root of $T$ is denoted as $\circ$. The set of internal nodes of $T$ is denoted as $I$, whereas its set of leaves is denoted as $L$. By definition, $\circ\in I$ and $I$ and $L$ partition $V$. From the definition of ORB-tree it also follows that $|L|=|I|=n$, hence $|V|=2n$. Note that $|E|=|V|-1$ because $T$ is a tree. We define $|T|:=|V|$. We use this later notation when we want to emphasize a direct relationship with the ORB-tree.

For $i,j\in V$, a path of length $l$ between $i$ and $j$ is a sequence $v_0,\ldots,v_l\in V$ such that $v_0=i$, $v_l=j$, and $\{v_k,v_{k+1}\}\in E$ for $0\le k<l$. Unless otherwise stated, we write $[i,j]$ to denote the set of edges in the shortest path in $T$ between $i$ and $j$. This path is unique because $T$ has no cycles. The depth of $i$, denoted $\text{depth}(i)$, is defined as $|[i,\circ]|$ i.e. the number of edges that connect $i$ with the root. We say that $i$ is an ancestor of $j$, or alternatively $j$ is a descendent of $i$, when $i\in[\circ,j]$. In particular, every node is an ancestor and a descendant from itself. Further, $(i\wedge j)$ denotes the so-called least-common ancestor to $i$ and $j$. This is the $v\in V$ that minimizes $|[v,\circ]|$, among all the nodes that are ancestors to both $i$ and $j$. 

We define
\[\ell(i,j):=\sum_{e\in[i,j]}\ell(e).\]
In addition, for $J\subset L$ and $i\in V$, define $\ell(J,i)$ as the column vector of dimension $|J|$ with entries $\ell(j,i)$, for $j\in J$.  $\ell(i,J)$ is the transpose of $\ell(J,i)$.

For each $i\in V$, $T(i)$ denotes the subtree of $T$ rooted at $i$. In particular, the vertex set of $T(i)$ is the subset of nodes in $T$ that descend from $i$, and its edge set is the subset of edges that connect two descendants of $i$. $L(i)$ denotes the leaf set of $T(i)$. Likewise, for each $e=\{i,j\}\in E$, if $i$ is closest to the root than $j$, $T(e)$ and $L(e)$ denote $T(i)$ and $L(i)$, respectively.

\section{Haar-like basis of ORB-trees}
\label{sec:Haarlike}

In this section we specialize the concept of Haar-like basis given in~\cite{gavnad10} to our setting of ORB-trees. The key new result in this section shows that the Haar-like basis of an ORB-tree interacts nicely with its covariance matrix (i.e. the strictly ultrametric matrix associated with the tree). This is somewhat unexpected because the covariance matrix is determined by the topology of the tree and its branch length function, whereas the basis is solely determined by the tree's topology.

To construct the Haar-like wavelets, it is convenient to represent the nodes in $I\setminus\{o\}$ momentarily as binary strings. With this convention, the (only) child of the root is $\varepsilon$---the so-called empty string. Further, the children of each node $v\in I\setminus\{o\}$ are $v0$ (i.e. the string $v$ with the character zero appended at the end) and $v1$ (i.e. $v$ with the character one appended at the end). 

\vspace{12pt}
\begin{tcolorbox}
\begin{definition}[Specialization from
\cite{gavnad10}]
The Haar-like basis associated with $T$ is the set of transformations $\{\varphi_v\}_{v\in I}$ defined as follows:
\[\varphi_o(i):=\frac{1}{\sqrt{|L|}},\text{ for all }i\in L;\]
and, for each $v\in I$ with $v\ne o$:
\[\varphi_v(i):=\begin{cases}
+\,\sqrt{\frac{|L(v1)|}{|L(v0)|\cdot|L(v)|}} &, i\in L(v0);\vspace{4pt}\\ 
-\,\sqrt{\frac{|L(v0)|}{|L(v1)|\cdot|L(v)|}} &, i\in L(v1);\\
0 &, \text{otherwise}.
\end{cases}\]
The Haar-like matrix associated with $T$ is the matrix $\Phi$ with columns $\varphi_v$, $v\in I$.
\label{def:pseudohaar}
\end{definition}
\end{tcolorbox}

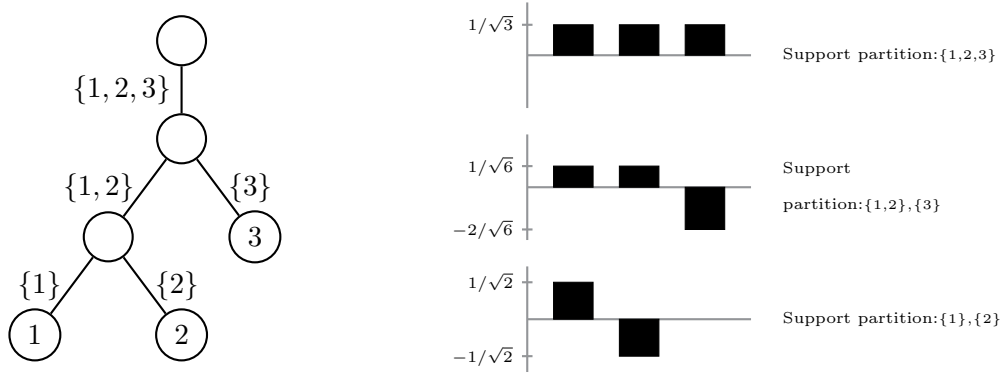
\begin{figure}
\centering
\begin{minipage}{0.49\textwidth}
\centering
\begin{tikzpicture}[scale=.65, node distance={15mm}, thick, main/.style = {draw, circle}] 
    \node[shape=circle,draw=black,minimum size=.65cm] (A) at (0,0) {};
    \node[shape=circle,draw=black,minimum size=.65cm] (B) at (0,-2) {};
    \node[shape=circle,draw=black] (C) at (1.5,-4) {3};
    \node[shape=circle,draw=black,minimum size=.65cm] (D) at (-1.5,-4) {};
    \node[shape=circle,draw=black] (E) at (-3,-6) {1};
    \node[shape=circle,draw=black] (F) at (0,-6) {2};

    \path [-] (A) edge node[left] {$\{1,2,3\}$} (B);
    \path [-] (B) edge node[left] {$\{1,2\}$} (D);
    \path [-] (B) edge node[right] {$\{3\}$} (C);
    \path [-] (D) edge node[left] {$\{1\}$} (E);
    \path [-] (D) edge node[right] {$\{2\}$} (F);
\end{tikzpicture} 
\end{minipage}
\begin{minipage}{0.49\textwidth}
\centering
\begin{tikzpicture}[scale=.7]
\draw[gray, thick] (0,0) -> (4.25,0);
\draw[gray, thick] (0,-1) -> (0,1);
\draw[gray, thick] (-.1,.577) -> (.1,.577);
\node[] at (-.7,.577) {$\scriptscriptstyle 1/\sqrt{3}$};
\node[text width=3cm] at (7,0) {\tiny Support partition:$\scriptscriptstyle\{1,2,3\}$};

\draw[gray, thick] (0,-2.5) -> (4.25,-2.5);
\draw[gray, thick] (0,-3.5) -> (0,-1.5);
\draw[gray, thick] (-.1,-2.1) -> (.1,-2.1);
\node[] at (-.7,-2.1) {$\scriptscriptstyle 1/\sqrt{6}$};
\draw[gray, thick] (-.1,-3.3) -> (.1,-3.3);
\node[] at (-.85,-3.3) {$\scriptscriptstyle -2/\sqrt{6}$};
\node[text width=3cm] at (7,-2.5) {\tiny Support partition:$\scriptscriptstyle \{1,2\},\{3\}$};

\draw[gray, thick] (0,-5) -> (4.25,-5);
\draw[gray, thick] (0,-6) -> (0,-4);
\draw[gray, thick] (-.1,-4.3) -> (.1,-4.3);
\node[] at (-.7,-4.3) {$\scriptscriptstyle 1/\sqrt{2}$};
\draw[gray, thick] (-.1,-5.7) -> (.1,-5.7);
\node[] at (-.85,-5.7) {$\scriptscriptstyle -1/\sqrt{2}$};
\node[text width=3cm] at (7,-5) {\tiny Support partition:$\scriptscriptstyle\{1\},\{2\}$};

\draw [fill=black] (.5,0) rectangle (1.25,.577);

\draw [fill=black] (1.75,0) rectangle (2.5,.577);

\draw [fill=black] (3,0) rectangle (3.75,.577);

\draw [fill=black] (.5,-2.5) rectangle (1.25,-2.1);

\draw [fill=black] (1.75,-2.5) rectangle (2.5,-2.1);

\draw [fill=black] (3,-2.5) rectangle (3.75,-3.3);

\draw [fill=black] (.5,-5) rectangle (1.25,-4.3);

\draw [fill=black] (1.75,-5) rectangle (2.5,-5.7);
\end{tikzpicture}
\end{minipage}
\caption{\textbf{Visualization of the Haar-like Wavelet basis associated with an ORB-tree.} Left: ORB-tree with leaves $1,2,3$. Edges are labeled by the subsets of leaves that descend from them. Right: Haar-like basis associated with the ORB-tree on the left.}
\label{fig:WaveletTree}
\end{figure}

To fix ideas see Figure~\ref{fig:WaveletTree}.

The terminology of basis in Definition~\ref{def:pseudohaar} is justified by the fact that
\begin{equation}
\text{if $u,v\in I$ then $\langle \varphi_u,\varphi_v\rangle=\indicator[u=v]$.} 
\label{ide:orthonormality}
\end{equation}
In particular, $\{\varphi_v\}_{v\in I}$ is an orthonormal basis of $\mathbb{R}^{|L|}$. (See the Appendix for a self-contained justification of the orthonormality of the Haar-like basis.) Note that the Haar-like matrix $\Phi$ has its rows indexed by $L$ and its columns indexed by $I$. Hence, since $|L|=|I|$, $\Phi$ is a square matrix, an orthonormal one.

Clearly, for each $v\in I$, $\varphi_v$ has $L(v)$ as its support. Further, because $\varphi_o$ is constant, the orthogonality property implies for $v\ne o$ that $\sum_{i\in L}\varphi_v(i)=0$. These two properties are essential for our arguments onwards.

The following definition is useful to understand the relationship between the Haar-like basis of an ORB-tree and its associated covariance matrix.

\vspace{12pt}
\begin{tcolorbox}
\begin{definition}
The trace branch length of $T$ is the function $\ell^*:E\to[0,\infty)$ defined as $\ell^*(e):=|L(e)|\,\ell(e)$, for each $e\in E$. 
\label{def:tracelength}
\end{definition}
\end{tcolorbox}

\vspace{12pt}
\begin{tcolorbox}
\begin{theorem}
If $v\in I$ then
$S\,\varphi_v=\diag(\ell^*(L,v))\,\varphi_v.$
\label{thm:wow!}
\end{theorem}
\end{tcolorbox}

\begin{proof}
Consider $v\in I$ and $j\in L$. If $j\notin L(v)$ then $(i\wedge j)=(v\wedge j)$, hence
\[(S\,\varphi_v)(j)=\sum_{i\in L(v)}\ell(i\wedge j,\circ)\,\varphi_v(i)=\ell(v\wedge j,\circ)\sum_{i\in L(v)}\varphi_v(i).\]
But, if $v=\circ$ then $\ell(v\wedge j,\circ)=0$. Instead, if $v\ne\circ$ then $\sum_{i\in L(v)}\varphi_v(i)=0$. In either case: $(S\,\varphi_v)(j)=0$. This shows the lemma for $j\notin L(v)$ because the entry associated with $j$ in  $\diag(\ell^*(L,v))\,\varphi_v$ is $\ell^*(v,j)\cdot\varphi_v(j)$, and the support of $\varphi_v$ is $L(v)$.

Next suppose that $j\in L(v)$. Then
\begin{align*}
(S\,\varphi_v)(j)
&=\sum_{i\in L(v)}\sum_{e\in[i\wedge j,o]}\ell(e)\,\varphi_v(i)\\
&=\sum_{i\in L(v)}\sum_{e\in[i\wedge
j,v]}\ell(e)\,\varphi_v(i)+\sum_{i\in L(v)}\varphi_v(i)\cdot\sum_{e\in[v,o]}\ell(e)\\
&=\sum_{i\in L(v)}\varphi_v(i)\sum_{e\in[i\wedge
j,v]}\ell(e),
\end{align*}
where for the last identity we have used that $\sum_{i\in L(v)}\varphi_v(i)=0$ if $v\ne\circ$, and $\sum_{e\in[v,o]}\ell(e)=0$ if $v=\circ$.  But note that if $i\in L(v)$ is such that $(i\wedge j)=v$ then $\sum_{e\in[i\wedge
j,v]}\ell(e)=0$. Instead, if $(i\wedge j)\ne v$ then $\varphi_v(i)=\varphi_v(j)$. As a result
\begin{align*}
(S\,\varphi_v)(j)
&=\varphi_v(j)\sum_{i\in L(v):\,i\wedge j\ne v}\,\sum_{e\in[i\wedge
j,v]}\ell(e)\\
&=\varphi_v(j)\sum_{e\in[j,v]}\,\,\sum_{i\in L(v):\,e\in[i\wedge j,v]}\ell(e)\\
&=\varphi_v(j)\sum_{e\in[j,v]}\ell(e)\,|L(e)|\\
&=\varphi_v(j)\,\ell^*(j,v),
\end{align*}
which shows the result.
\end{proof}

It follows from the theorem that for each $u,v\in I$:
\begin{equation}
(\Phi'S\Phi)(u,v)=\varphi_u'S\varphi_v=\sum_{i\in L(v)\cap L(u)}\varphi_u(i)\,\varphi_v(i)\,\ell^*(i,v).
\label{ide:varphiu'Svarphiv}
\end{equation}
The importance of the diagonal of $\Phi'S\Phi$ in the discussion ahead, motivates to define for $v\in I$ the quantities
\begin{equation}
\lambda_v:=(\Phi'S\Phi)(v,v)=\sum_{i\in L(v)}\varphi_v^2(i)\,\ell^*(i,v).
\label{def:lambdav}
\end{equation}

For $v\in I$, because $\varphi_v$ has $L(v)$ as its support and $\|\varphi_v\|_2=1$, $\lambda_v$ is a weighted average of the trace branch length between each leaf in $L(v)$ and $v$. In particular, since $L(u)\supset L(v)$ when $u$ is an ancestor of $v$, the closer the internal node $v$ is to the root, the more terms are averaged. (This emulates the averaging at different scales that the standard Haar wavelet transform does to a continuous signal.) Furthermore, since $\ell^*(e)=\ell(e)>0$ when $e$ joins a leave with its parent, $\lambda_v>0$.

On the other hand, because $(\Phi'S\Phi)(u,v)=0$ when $u,v\in I$ are such that $L(u)\cap L(v)=\emptyset$, the identity in (\ref{ide:varphiu'Svarphiv}) suggests that the Haar-like matrix can be used to sparsify the covariance matrix of the ORB-tree. The following result is critical to assess how effective this sparsification is in practice.

\vspace{12pt}
\begin{tcolorbox}
\begin{lemma}
For all $u,v\in V$, $L(u)\cap L(v)\ne\emptyset$ if and only if $u$ is an ancestor of $v$ or vice versa.
\label{lem:LuLvnotempty}
\end{lemma}
\end{tcolorbox}

\begin{proof} 
If $u$ is an ancestor of $v$ then $L(v)\subset L(u)$; in particular, $L(u)\cap L(v)=L(v)\ne\emptyset$. The same conclusion applies if $v$ is an ancestor of $u$. Conversely, suppose that $L(u)\cap L(v)\ne\emptyset$. Without loss of generality assume that $u\ne v$. From the hypothesis, there is $w\in L$ that descends from both $u$ and $v$. But, since there is a unique path from $w$ to $\circ$, $u$ and $v$ must be both in this path; in particular, either $u$ is an ancestor of $v$ or vice versa.
\end{proof}

\section{Sparisfication of Covariance Matrices of ORB-trees}
\label{sec:Sparisfication}

In this section we quantify how much of the covariance matrix of an ORB-tree can be sparsified by its Haar-like matrix. To state our main result we require the following definitions.

\begin{tcolorbox}
\begin{definition}
The average subtree size of $T$ is the quantity, $\avg(T):=\frac{1}{|T|}\sum\limits_{v\in V}|T(v)|$.
\end{definition}
\end{tcolorbox}

\begin{tcolorbox}
\begin{definition}
The internal and external path lengths of $T$ are the quantities defined as $\IPL(T):=\sum_{v\in I}\text{depth}(v)$ and $\EPL(T):=\sum_{v\in L}\text{depth}(v)$, respectively~\cite{sedfla13}. The total path length of $T$ is the quantity $\TPL(T):=\IPL(T)+\EPL(T)$.
\end{definition}
\end{tcolorbox}

We note the relationship:
\begin{equation}
\avg(T)=1+\frac{\TPL(T)}{|T|},
\label{ide:avgTPL}
\end{equation}
because
\[\TPL(T)
=\sum_{v\in V}\sum_{u\in V\setminus\{\circ\}}\indicator[v\in T(u)]
=\sum_{u\in V\setminus\{\circ\}}|T(u)|
=\left\{\sum_{u\in V}|T(u)|\right\}-|T|.
\]

\begin{tcolorbox}
\begin{definition}
The interior of $T$ is the tree $\mathring{T}$ obtained by trimming the leaves of $T$.
\end{definition}
\end{tcolorbox}

Clearly, $\IPL(T)=\TPL(\mathring{T})$.

As mentioned earlier, the identity in (\ref{ide:varphiu'Svarphiv}) guarantees that some entries of $\Phi'S\Phi $ vanish. The following result estimates the least number of such entries. Our lower bound is independent of the branch lengths and depends---only---on the tree topology.

\vspace{12pt}
\begin{tcolorbox}
\begin{theorem}
If $\zeta$ denotes the fraction of vanishing entries in $\Phi'S\Phi$ then
\[\zeta\ge1+\frac{1}{|L|}-2\frac{\avg(\mathring{T})}{|\mathring{T}|}=1-\frac{1}{|L|}-2\frac{\TPL(\mathring{T})}{|\mathring{T}|^2}.\] 
\label{thm:zetaexplicit}
\end{theorem}
\end{tcolorbox}

\begin{proof}
Recall that $|I|=|L|$ and, for $v\in I$, the support of $\varphi_v$ is $L(v)$. Hence, from the identity in~(\ref{ide:varphiu'Svarphiv}), $(\Phi'S\Phi)(u,v)=0$ when $u,v\in I$ and $L(u)\cap L(v)=\emptyset$. As a result, using that $L(u)\ne\emptyset$ when $u\in I$, and Lemma~\ref{lem:LuLvnotempty}, we obtain that
\begin{align*}
|I|^2\zeta
&\ge |I|^2-|\{(u,v)\in I\times I\text{ such that }L(u)\cap L(v)\ne\emptyset\}|\\
&= |I|^2-|I|\\
&\qquad-2|\{(u,v)\in I\times I\text{ such that $v\ne u$ descends from $u$ and }L(u)\cap L(v)\ne\emptyset\}|\\
&= |I|^2-|I|-2\sum_{u\in I}\big(|\mathring{T}(u)|-1\big)\\
&= |I|^2+|I|-2\sum_{u\in I}|\mathring{T}(u)|,
\end{align*}
Since $|I|=|L|=|\mathring{T}|=|T|/2$, $|L|^2\zeta\ge|L|^2+|L|-|\mathring{T}|\cdot\avg(\mathring{T})$, which shows the inequality in the theorem. The alternative lower-bound for $\zeta$ follows by applying the identity in equation~(\ref{ide:avgTPL}) to $\mathring{T}$, completing the proof of the theorem.
\end{proof}

It follows from the first lemma in~\cite[Section 6.4]{sedfla13} that for an ORB-Tree T, $\EPL(T)-\IPL(T)=2|I|-1$, which together with the previous theorem let us conclude the following asymptotic result.

\begin{tcolorbox}
\begin{corollary}
If either $\avg(\mathring{T})\ll|\mathring{T}|$, $\TPL(\mathring{T})\ll|\mathring{T}|^2$, $\IPL(T)\ll|I|^2$, or $\EPL(T)\ll|L|^2$ as $|T|\to\infty$, then $\zeta=1-o(1)$.
\label{cor:zeta21}
\end{corollary}
\end{tcolorbox}

In other words, if $T$ grows so that either of the asymptotic inequalities in the above corollary applies, then an asymptotically negligible fraction of the off-diagonal entries in $\Phi'S\Phi$ will be non-zero. 

The last asymptotic condition in the corollary (i.e., that $\EPL(T)\ll|L|^2$) is of relevance in phylogenetic studies. In that context, the external path length of a tree is called its Sackin's index~\cite{BluOli05,CorArnRosRot20,KinRos21}. This index is used as a measure of the imbalance of phylogenetic trees. In particular, since phylogenetic trees are neither too balanced nor too imbalanced~\cite{Ald01,BluFra06}; Haar-like bases should be rather effective in sparsifying covariance matrices of phylogenetic trees in practice. We come back to this point in Section~\ref{sec:treeoflife}.

\subsection{Sparsification of covariance matrices of maximally balanced ORB-trees}

In our context, the following definition gives the most balanced topology among the ORB-trees.

\begin{tcolorbox}
\begin{definition}[Perfect Binary Trees]
A perfect binary tree is an ORB-tree in which all leaves have the same depth.
\label{def:perfectbin}
\end{definition}
\end{tcolorbox}

To fix ideas see Figure~\ref{fig:bin4}. 

Let $T$ be a perfect binary tree of height $(h+1)$. In particular, {$|L|=2^{h}$} and {$|V|=2^{h+1}$.} At level $k\ge1$, $T$ contains $2^{k-1}$ nodes, each of which is the root of a perfect binary tree of height $(h-k)$. Since a perfect binary tree of height $h$ contains $(2^{h+1}-1)$ nodes, and the interior of a perfect binary tree of height $(h+1)$ is a perfect binary tree of height $h$:
\[\avg(\mathring{T})=\frac{2^h+\sum\limits_{k=1}^h 2^{k-1}\cdot (2^{h-k+1}-1)}{2\cdot2^{h-1}}=h+2^{-h}\ll |\mathring{T}|.\]
\color{black}
In particular, due to Theorem~\ref{cor:zeta21}, we can conclude that the Haar-like matrix can be used to asymptotically annihilate (via a similarity transformation) the off-diagonal entries of the covariance matrix of a perfect binary tree  as its height tends to infinity.

\subsection{Sparsification of covariance matrices of maximally imbalanced ORB-trees}

The following definition provides what we may regard as the most unbalanced topology among the ORB-trees.

\begin{tcolorbox}
\begin{definition}[Binary Caterpillar Trees]
The binary caterpillar tree of height $h\ge1$ is the ORB-tree with vertices $\circ,1,\ldots,h,1',\dots,(h-1)'$ and edges of the form $\{\circ,1\}$, $\{i,i+1\}$, for $i=1,\ldots,(h-1)$, and $\{j,j'\}$ for $j=1,\ldots,(h-1)$. 
\label{def:bincat}
\end{definition}
\end{tcolorbox}

To fix ideas see Figure~\ref{fig:cateh4}.

Let $T$ be a binary Caterpillar tree of height $h$. In particular, $|V|=2h$, $|L|=h$, and each node in level $k\ge1$ is the root of a snake binary subtree of height $(h-k)$. Therefore, each internal node has as children one leaf node and one internal node that is the root of a binary caterpillar subtree, and the subgraph of internal nodes is a path. As a result:
\[\avg(\mathring{T})=\frac{2\sum\limits_{k=1}^{h-1}\big(h-k\big)}{2(h-1)}=\frac{h}{2}+1\sim\frac{|\mathring{T}|}{2}.\]
\color{black}
Hence, the lower-bound provided by Theorem~\ref{thm:zetaexplicit} is trivial, and we cannot guarantee that the Haar-like matrix associated with a sizeable binary caterpillar tree annihilates its off-diagonal entries in any significant way.

\subsection{Sparsification of covariance matrices of large random ORB-trees}

Perfect binary trees and caterpillar trees are opposite extremes of how balanced (or imbalanced) ORB-trees can be. It is therefore unclear how much sparsification the Haar-like matrix of a large but generic ORB-tree can induce on its covariance matrix. To address this issue we consider a natural ensemble of random ORB-trees

In what follows, $\TT$ denotes a uniformly at random ORB-tree with $|I|$ internal nodes. Such trees may be generated using the Catalan distribution~\cite[Section 6.7]{sedfla13}. This probability model produces full binary trees with a given number of internal nodes; which we may turn into an ORB-tree by appending their root to a new one. 

Let $\mathbb{S}$ denote the covariance matrix of $\TT$, and $\zeta$ the number of zeroes in the random matrix $\Phi'\mathbb{S}\Phi$, where $\Phi$ is the Haar-like matrix associated with $\TT$. It turns out that the mean and variance of the internal path length of $\TT$ are given by 
\begin{align}
\EE\big(\IPL(\TT)\big)
&\sim \sqrt{\pi}|I|^{3/2};
\label{ide:EEIPLTT}\\
\VV\big(\IPL(\TT)\big)
&\sim\left(\frac{10}{3}-\pi\right)|I|^3.
\label{ide:VVIPLTT}
\end{align}
The identity in equation~(\ref{ide:EEIPLTT}) follows from~\cite[Proposition VII.3.]{flased13}. The identity in~(\ref{ide:VVIPLTT}) may be regarded a refinement of~\cite[Note VII.12]{flased13}. 

As the following result implies, the Haar-like basis of most large ORB-trees should be highly effective in sparsifying their covariance matrix.

\begin{tcolorbox}
\begin{corollary}
If $\TT$ is a uniformly at random ORB-tree with $|I|$ internal nodes then $\zeta=1-o(1)$ with overwhelmingly high probability, as $|I|\to\infty$.
\label{cor:randomT}
\end{corollary}
\end{tcolorbox}

\begin{proof}
Let $t>0$. Let $\mu$ and $\sigma^2$ denote the mean and variance of $\IPL(\TT)$, respectively. Due to Cantelli's inequality (a one sided version of the well-known Chebyshev's inequality): $\PP\big(\IPL(\TT)\ge\mu+t\sigma\big)\le(1+t^2)^{-1}$. But $(\mu+t\sigma)=\Omega\big(t|I|^{3/2}\big)$ because of equations~(\ref{ide:EEIPLTT})-(\ref{ide:VVIPLTT}). In particular, there is a constant $c>0$ such that
\[\PP\left(\frac{\IPL(\TT)}{|I|^2}\le\frac{ct}{\sqrt{|I|}}\right)\ge\frac{t^2}{1+t^2}.\]
So, if $t\to\infty$ so that $t=o(\sqrt{|I|})$ then $\frac{\IPL(\TT)}{|I|^2}=o(1)$ with a probability converging to one as $|I|\to\infty$. The result now follows from Corollary~\ref{cor:zeta21}.
\end{proof}
\color{black}

\subsection{Fast Sparsification Algorithm} 

A non-trivial challenge to storing and manipulating large strictly ultrametric matrices is that they are almost always fully dense in applications. Further, in the context of phylogenetic covariance matrices, the ORB-trees associated with such matrices are formed in advance. This allows us to sparsify these matrices without computing them, or even storing them in computer memory. It also allows us to anticipate which entries may remain nonzero after sparsification. In fact, due to equation~(\ref{ide:varphiu'Svarphiv}) and Lemma~\ref{lem:LuLvnotempty}, all that is required to sparsify these matrices from their ORB-tree is to precompute the leaves that descend from each internal node (i.e., the sets $L(v)$, with $v\in I$) and the trace branch length between them (Definition~\ref{def:tracelength}). This can be achieved with two postorder traversals of the ORB-tree. We convey these ideas in the following pseudo-code (Algorithm~\ref{alg:sparsify}), which is fully coded and available on \href{https://github.com/edgor17/Sparsify-Ultrametric}{GitHub}.


\begin{algorithm}
\begin{small}
\caption{Phylogenetic covariance matrix sparsification}
\label{alg:sparsify}
\begin{algorithmic}
\STATE{\textbf{Input.} ORB-tree $T$ with covariance matrix $S$}
\STATE{\textbf{Output.} Only possibly non-zero entries in $\Phi'S\Phi$}
\FOR{$v\in I$ in postorder traversal of $T$}
\FOR{$i\in L$}
\IF{$v=\circ$}
\STATE{$\varphi_o(i)\gets\frac{1}{\sqrt{|L|}}$}
\ELSIF{$i\in L(v0)$} \STATE{$\varphi_v(i)\gets+\,\sqrt{\frac{|L(v1)|}{|L(v0)|\cdot|L(v)|}}$}
\STATE{$\ell^*(v,i)\gets\ell^*(i,v0)+|L(v0)|\cdot\ell(v0,v)$}
\ELSIF{$i\in L(v1)$}
\STATE{$\varphi_v(i)\gets-\,\sqrt{\frac{|L(v0)|}{|L(v1)|\cdot|L(v)|}}$}
\STATE{$\ell^*(v,i)\gets\ell^*(i,v1)+|L(v1)|\cdot\ell(v1,v)$}
\ELSE
\STATE{$\varphi_v(i)\gets0$}
\ENDIF
\ENDFOR
\ENDFOR
\FOR{$v\in I$ in postorder traversal of $T$}
\WHILE{$\text{parent}(v)\neq \emptyset$}
\STATE{$u\gets\text{parent}(v)$}
\STATE{$M(u,v)\gets\sum\limits_{i\in L(v)\cap L(u)}\varphi_u(i)\,\varphi_v(i)\,\ell^*(v,i)$}
\ENDWHILE
\ENDFOR
\RETURN{$M(u,v)$ for $u,v\in I$ such that $L(u)\cap L(v)\ne\emptyset$}
\end{algorithmic}
\end{small}
\end{algorithm}

\section{Spectrum of Covariance Matrices of Trace-balanced ORB-trees}
\label{sec:Spectrum}

While Theorem \ref{thm:zetaexplicit} guarantees that some entries in $\Phi'S\Phi$ vanish---regardless of branch lengths, additional constraints on the latter can lead to further sparsification. The next definition identifies the class of ORB-trees whose Haar-like basis fully sparsifies (i.e, diagonalizes) their associated strictly ultrametric matrix.

\vspace{12pt}
\begin{tcolorbox}
\begin{definition}
$T$ is called trace-balanced at a node $v$ when, for all $i,j\in L(v)$, $\ell^*(i,v)=\ell^*(j,v)$. $T$ is called trace-balanced when it is trace-balanced at each $v\in I\setminus\{\circ\}$.
\label{def:tracebalance}
\end{definition}
\end{tcolorbox}

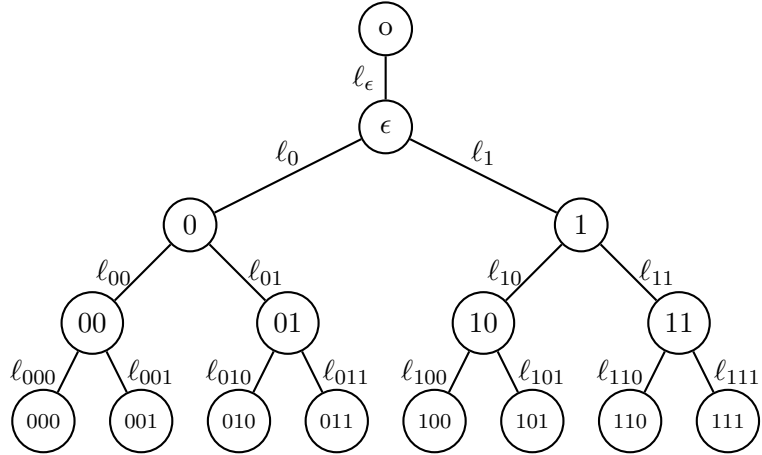
\begin{figure}
\begin{center}
\begin{tikzpicture}[scale=.65, node distance={15mm}, thick, main/.style = {draw, circle}] 

\node[shape=circle,draw=black,minimum size=.7cm] (A) at (-1,0) {o} ;
 \node[shape=circle,draw=black,minimum size=.7cm] (B) at (-1,-2) {$\epsilon$} ;
 
\node[shape=circle,draw=black,minimum size=.7cm] (C) at (-5,-4) {0} ;
\node[shape=circle,draw=black,minimum size=.7cm] (D) at (3,-4) {1};    

\node[shape=circle,draw=black] (E) at (-7,-6) {00} ;
 \node[shape=circle,draw=black] (F) at (-3,-6) {01} ;
\node[shape=circle,draw=black] (G) at (1,-6) {10};
\node[shape=circle,draw=black] (H) at (5,-6) {11};

\node[shape=circle,draw=black] (I) at (-8,-8) {$\scriptstyle 000$} ;
 \node[shape=circle,draw=black] (J) at (-6,-8) {$\scriptstyle 001$} ;
\node[shape=circle,draw=black] (K) at (-4,-8) {$\scriptstyle 010$};
\node[shape=circle,draw=black] (L) at (-2,-8) {$\scriptstyle 011$};
\node[shape=circle,draw=black] (M) at (0,-8) {$\scriptstyle 100$};
\node[shape=circle,draw=black] (N) at (2,-8) {$\scriptstyle 101$};
\node[shape=circle,draw=black] (O) at (4,-8) {$\scriptstyle 110$};
\node[shape=circle,draw=black] (P) at (6,-8) {$\scriptstyle 111$};

    \path [-] (A) edge node[left] {$\ell_{\epsilon}$} (B);
    \path [-] (B) edge node[above] {$\ell_1$} (D);
    \path [-] (B) edge node[above] {$\ell_0$} (C);
        \path [-] (C) edge node[left] {$\ell_{00}$} (E);
        \path [-] (C) edge node[right] {$\ell_{01}$} (F);
        \path [-] (D) edge node[left] {$\ell_{10}$} (G);
        \path [-] (D) edge node[right] {$\ell_{11}$} (H);
        \path [-] (E) edge node[left] {$\ell_{000}$} (I);
        \path [-] (E) edge node[right] {$\ell_{001}$} (J);
       \path [-] (F) edge node[left] {$\ell_{010}$} (K);
        \path [-] (F) edge node[right] {$\ell_{011}$} (L);
       \path [-] (G) edge node[left] {$\ell_{100}$} (M);
        \path [-] (G) edge node[right] {$\ell_{101}$} (N);
               \path [-] (H) edge node[left] {$\ell_{110}$} (O);
        \path [-] (H) edge node[right] {$\ell_{111}$} (P);
\end{tikzpicture}
\caption{\textbf{Visualization of a perfect binary tree of height 4.} Such tree is trace-balanced if and only if $\ell_\alpha=\ell_\beta$ for each pair of binary strings $\alpha$ and $\beta$ of the same length. If all these lengths are strictly positive, the eigenvalues of its covariance matrix are $\ell_{000}$ (multiplicity 4), $\ell_{000}+2\,\ell_{00}$ (multiplicity 2), $\ell_{000}+2\,\ell_{00}+4\,\ell_0$ (multiplicity 1), and $\ell_{000}+2\,\ell_{00}+4\,\ell_0+8\,\ell_\epsilon$ (multiplicity 1). Otherwise, some multiplicities need to be added up.}
\label{fig:bin4}
\end{center}
\end{figure}

\begin{figure}
\begin{center}
\begin{tikzpicture}[scale=.55, node distance={15mm}, thick, main/.style = {draw, circle}] 
\node[shape=circle,draw=black,minimum size=.6cm] (A) at (-1,0) {o} ;
 \node[shape=circle,draw=black,minimum size=.7cm] (B) at (-1,-2) {$1$} ;
 
\node[shape=circle,draw=black,minimum size=.6cm] (C) at (-3,-4) {$1'$} ;
\node[shape=circle,draw=black,minimum size=.7cm] (D) at (1,-4) {2};    
      
\node[shape=circle,draw=black,minimum size=.6cm] (G) at (-1,-6) {2'};
\node[shape=circle,draw=black,minimum size=.7cm] (H) at (3,-6) {$3$};

\node[shape=circle,draw=black,minimum size=.6cm] (O) at (1,-8) {$3'$};
\node[shape=circle,draw=black,minimum size=.6cm] (P) at (5,-8) {$4$};

    \path [-] (A) edge node[left] {$\ell_0$} (B);
    \path [-] (B) edge node[above right] {$\ell_1$} (D);
    \path [-] (B) edge node[above left] {$\ell_1'$} (C);
        \path [-] (D) edge node[above left] {$\ell_{2}'$} (G);
        \path [-] (D) edge node[above right] {$\ell_{2}$} (H);
                      \path [-] (H) edge node[above left] {$\ell_{3}'$} (O);
        \path [-] (H) edge node[above right] {$\ell_{3}$} (P);
\end{tikzpicture}
\caption{\textbf{Visualization of a binary caterpillar tree of height 4.}  This tree is trace-balanced if and only if $\ell_3'=\ell_3$, $\ell_2'=\ell_3+2\,\ell_2$, and $\ell_1'=\ell_3+2\,\ell_2+3\,\ell_1$. In such case, if $\ell_3,\ell_2,\ell_1,\ell_0>0$ then its covariance matrix spectrum is $\{\ell_3, \ell_3+2\,\ell_2, \ell_3+2\,\ell_2+3\,\ell_1, \ell_3+2\,\ell_2+3\,\ell_1+4\,\ell_0\}$, and each eigenvalue is simple.}
\label{fig:cateh4}
\end{center}
\end{figure}
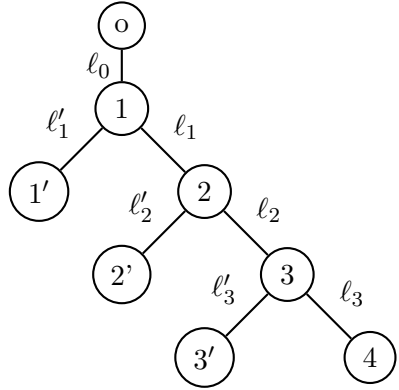

We note that a tree is always trace-balanced at a leaf. Also, if an ORB-tree is trace-balanced at the child of its root then it is also trace-balanced at the root. This is why the definition of trance-balanced trees only considers nodes in $I\setminus\{o\}$. See Figures~\ref{fig:bin4}-\ref{fig:cateh4} for depictions of trace-balanced trees.

The following two results show the relevance of the above definition in terms of the eigenvalues of the ultrametric matrix associated with an ORB-tree. 

\vspace{12pt}
\begin{tcolorbox}
\begin{lemma}
If $v\in I$ then $\varphi_v$ is an eigenvector of $S$ if and only if $T$ is trace-balanced at $v$, in which case the eigenvalue associated with $\varphi_v$ is $\ell^*(i,v)$, for any $i\in L(v)$.
\label{lem:varphiveig}
\end{lemma}
\end{tcolorbox}

\begin{proof}
Fix $v\in I$. Due to Theorem~\ref{thm:wow!}, $(S\varphi_v)(i)=\ell^*(i,v)\,\varphi_v(i)$, for each $i\in L$. This shows the lemma because $\varphi_v(i)>0$ if and only if $i\in L(v)$.
\end{proof}

Since the covariance matrix $S$ has dimensions $|L|\times|L|$, but $|I|=|L|$ because $T$ is an ORB-tree, the following result is immediate from the previous lemma.

\vspace{12pt}
\begin{tcolorbox}
\begin{corollary}
The Haar-like basis of $T$ diagonalizes its covariance matrix if and only if $T$ is trace-balanced. In this case, the spectrum of $S$ is
\[\sigma(S)=\bigcup\limits_{v\in I}\{\ell^*(v,i)\text{ for any }i\in L(v)\},\]
and the multiplicity of $\ell^*(v,i)$ is $\big|\{u\in I:\,\ell^*(v,i)=\ell^*(u,j),\text{ for some }j\in L(u)\}\big|$.
\label{cor:varphiveig}
\end{corollary}
\end{tcolorbox}

Due to Corollary~\ref{cor:varphiveig}, we can assert the following.

\begin{tcolorbox}
\begin{corollary}
The perfect binary tree of height $h$ is trace-balanced if and only if it has constant branch lengths at each level. In this case, if $\ell_j$ denotes the common length of the edges that connect a node at depth $j$ with another at depth $(j+1)$, the spectrum of the associated covariance matrix is
\[\sigma(S)=\left\{\sum_{k=j}^{h-1}2^{h-k}\ell_k,\text{ with }j=0,\ldots,h-1\right\}.\]
Furthermore, the multiplicity of the eigenvalue $\lambda=\sum\limits_{k=j}^{h-1}2^{h-k}\ell_k$ is $\max\left\{1,\sum\limits_{j\in\Lambda}2^{j-1}\right\}$, where

\[\Lambda:=\left\{j'\in\{0,\ldots,h-1\}\text{ such that }\sum_{k=j'}^{h-1}2^{h-k}\ell_k=\sum_{k=j}^{h-1}2^{h-k}\ell_k\right\}.\]
\label{cor:fullconstbranchperlevel}
\end{corollary}
\end{tcolorbox}

Next, consider the binary caterpillar tree from Definition \ref{def:bincat}. In particular, its internal and leaf set are $I=\{\circ,1,\ldots,h-1\}$ and $L=\{1',\ldots,(h-1)',h\}$, respectively. Let $\ell_0$ denote the branch length of $\{\circ,1\}$, $\ell_i$ the length of $\{i,i+1\}$ for $i=1,\ldots,(h-1)$, and $\ell_j'$ the branch length of $\{j,j'\}$ for $j=1,\ldots,(h-1)$. Due to Corollary~\ref{cor:varphiveig}, we have the following result.

\begin{tcolorbox}
\begin{corollary}
The Caterpillar tree of height $h$ is trace-balanced if and only if $\ell_j'=\sum\limits_{k=j}^{h-1}(h-k)\cdot\ell_k$, for $j=1,\ldots,h-1$. In this case, the eigenvalues of the associated covariance matrix $S$ are as follows, repeated according to their multiplicity: $\ell_0'\ge\ell_1'\ge\cdots\ge\ell_{h-1}'$, where $\ell_0':=\sum\limits_{k=0}^{h-1}(h-k)\cdot\ell_k$.
\label{cor:balancedCaterpillar}
\end{corollary}
\end{tcolorbox}

We finish this section with a definition and result that characterizes the possible spectrums of covariance matrices of trace-balanced trees. Because the result's proof is constructive, it can be used to form strictly ultrametric matrices with the desired spectrum and multiplicities.

\begin{tcolorbox}
\begin{definition}
In a tree $T$, a function $f:I\to[0,\infty)$ is called decreasing when, for all distinct $u,v\in I$, if $u$ is an ancestor of $v$ then $f(u)\ge f(v)$. In addition, $f$ is called strictly positive at the fringe when $f(u)>0$ whenever $u$ is a parent of a leaf.
\end{definition}
\end{tcolorbox}

\vspace{12pt}
\begin{tcolorbox}
\begin{corollary}
In a trace-balanced ORB-tree $T$ the function $v\longrightarrow\ell^*(v,i)$, with $v\in I$ and any $i\in L(v)$, is decreasing, and strictly positive at the fringe. Conversely, given any ORB-tree topology $T$ and decreasing function $f:V\to[0,\infty)$ that is strictly positive at the fringe, there is a branch length function $\ell:E\to[0,\infty)$ such that $\sigma(S)=f(I)$. Furthermore, the multiplicity of $\lambda\in\sigma(S)$ is $|f^{-1}(\{\lambda\})|$.
\label{cor:decfct}
\end{corollary}
\end{tcolorbox}
\begin{proof}
From the definitions of trace-balanced and ORB-tree, it is immediate that the transformation $v\in I\longrightarrow \ell^*(v,i)$, with $i\in L(v)$, is well-defined and strictly positive at the fringe of $T$. Also, $f$ is decreasing because if $u$ is an ancestor of $v$ then, for each $i\in L(v)$: $\ell^*(u,i)=\ell^*(u,v)+\ell^*(v,i)\ge\ell^*(v,i)$. This shows the first statement in the corollary.

For the second statement consider an ORB-tree topology $T=(V,E)$ and function $f:I\rightarrow [0,\infty)$ that is both decreasing and strictly positive at the fringe. Due to Corollary~\ref{cor:varphiveig}, it suffices to show that there is a branch length function $\ell:E\to[0,+\infty)$ such that $f(v)=\ell^*(v,i)$, for all $i\in L(v)$. To do so, let $e=\{u,v\}\in E$ be so $\text{depth}(u)<\text{depth}(v)$. Define
\begin{equation}
\ell(e):=
\begin{cases}
f(u) &, v\in L;\\
\frac{f(u)-f(v)}{|L(e)|} &, v\in I.
\end{cases}
\end{equation}
Observe that if $v\in L$ then $|L(e)|=1$ so $f(u)=\ell^*(e)$. Further, $\ell(e)>0$ because $f$ is strictly positive at the fringe of $T$. Instead, if $v\in I$ then $f(u)=f(v)+\ell(e)\cdot|L(e)|=f(v)+\ell^*(e)$, and $\ell(e)\ge0$ because $f$ is decreasing. In particular, if we extend the domain of $f$ to all of $V$ defining $f(v):=0$ for $v\in L$ then, for all $e=\{u,v\}\in E$ such that $\text{depth}(u)<\text{depth}(v)$: $f(u)=f(v)+\ell^*(e)$. From this, a simple inductive argument on the difference $d:=\text{depth}(v)-\text{depth}(u)>0$ shows that $f(u)-f(v)=\ell^*(u,v)$; implying that $f(u)=\ell^*(u,i)$, for all $i\in L(u)$, as claimed.
\end{proof}

\section{Spectrum Approximation in Roughly Trace-balanced ORB-trees}
\label{sec:Approx}

We know that the Haar-like wavelet associated with an internal node of an ORB-tree is an eigenvector of its covariance matrix if and only if the node is trace-balanced (Lemma~\ref{lem:varphiveig}). On the other hand, the Haar-like matrix of an ORB-tree can sometimes sparsify its covariance matrix significantly (Theorem~\ref{thm:zetaexplicit}). Together, these two facts suggest that the diagonal entry in $\Phi'S\Phi$ associated with an ``approximately'' trace-balanced internal node should be near the spectrum of $S$.
In this section, we formalize this intuition quantifying what it is required for an internal node to be approximately trance-balanced.

In what follows, for a given function $x:L\to\RR$ and non-empty $J\subset L$, we define the average value and variance of $x$ over $J$ as the quantities
\begin{align*}
\avg(x;J) &:=\frac{1}{|J|}\sum_{j\in J}x(j);\\
\var(x;J) &:=\frac{1}{|J|}\sum_{j\in J}\big(x(j)-\avg(x;J)\big)^2.
\end{align*}

In addition, for each $v\in V$, let $\text{parent}(v)$ denote the parent of node $v$ in $T$. Define
\[\rho_v:=\frac{|L(v)|}{|L\big(\text{parent}(v)\big)|}.\]

The following result aids in formalizing the intuition mentioned earlier.

\begin{tcolorbox}
\begin{lemma}
If $A$ is a symmetric matrix of dimensions $n\times n$ then, for all $\lambda\in\RR$:
\[\text{distance}\big(\lambda,\sigma(A)\big)\le\min_{x\in\RR^n:\,\|x\|_2=1}\|(A-\lambda)x\|_2.\]
\label{lem:phew!}
\end{lemma}
\end{tcolorbox}

\begin{proof}
Since $A$ is symmetric, all its eigenvalues are real and $\RR^n$ has an orthonormal basis of eigenvectors $v_1,\ldots,v_n$. Say $Av_i=\lambda_iv_i$. In particular, if $\lambda\in R$ and $x\in\RR^n$ then
\[(A-\lambda)x = \sum_{i=1}^n(\lambda_i-\lambda)\langle x,v_i\rangle v_i.\]
So, if $\|x\|_2=1$ then
\begin{align*}
\|(A-\lambda)x\|_2
&= \sqrt{\sum_{i=1}^n(\lambda_i-\lambda)^2\langle x,v_i\rangle^2}\\
&\ge\min_{\lambda'\in\sigma(A)}|\lambda'-\lambda|\cdot\sqrt{\sum_{i=1}^n\langle x,v_i\rangle^2}=\text{distance}\big(\lambda,\sigma(A)\big),
\end{align*}
which shows the lemma.
\end{proof}

Next, we provide a sufficient condition for $\lambda_v$, with $v\in I$, to be a good approximation of an eigenvalue of $S$. We also quantify explicitly the cosine between $S\varphi_v$ and $\lambda_v\varphi_v$ to assess how close $\varphi_v$ is to be an eigenvector of $S$. 

In the following result we use the notation: $\neg 0=1$ and $\neg 1=0$.

\begin{tcolorbox}
\begin{theorem}
If $v\in I$ then $\lambda_v=\rho_{v1}\cdot\overline{\ell^*\big(L(v0),v\big)}+\rho_{v0}\cdot\overline{\ell^*\big(L(v1),v\big)}$, and
\[\text{distance}\big(\lambda_v,\sigma(S)\big)\qquad\qquad\qquad\qquad\qquad\qquad\qquad\qquad\qquad\qquad\qquad\qquad\qquad\qquad\quad\]
\vspace{-30pt}
\begin{align*}
\color{red}
&\le\|(S-\lambda_v)\varphi_v\|_2\\
&=\sqrt{\rho_{v0}\cdot\rho_{v1}\cdot\Big\{\overline{\ell^*(L(v1),v)}-\overline{\ell^*(L(v0),v)}\Big\}^2+\sum_{\alpha\in\{0,1\}}\rho_{v\alpha}\cdot\var\big(\ell^*(L,v);L(v\neg \alpha)\big)}.
\end{align*}
Furthermore,
\[\cos(S\varphi_v,\lambda_v\varphi_v) = \frac{1}{\sqrt{1+\left\{\frac{\|(S-\lambda_v)\varphi_v\|_2}{\lambda_v}\right\}^2}}.\]
\label{thm:eigest}
\end{theorem}
\end{tcolorbox}

\begin{proof}
Fix $v\in I$. To make the $\lambda_v$ more explicit, observe that if $x:L\to\RR$ is a function (or vector) then
\begin{equation}
\sum_{i\in L}\varphi_v^2(i)\cdot x(i)=\rho_{v1}\cdot\avg\big(x;L(v0)\big)+\rho_{v0}\cdot\avg\big(x;L(v1)\big).
\label{lem:avgOFavg} 
\end{equation}
In particular, due to Theorem~\ref{thm:wow!}:
\begin{equation}
\lambda_v=\varphi_v'\,S\varphi_v=\sum_{i\in L(v)}\varphi_v^2(i)\cdot\ell^*(v,i)=\rho_{v1}\cdot\overline{\ell^*\big(L(v0),v\big)}+\rho_{v0}\cdot\overline{\ell^*\big(L(v1),v\big)},
\label{ide:lambda*}
\end{equation}
which shows the first identity in the theorem.

On the other hand, Lemma~\ref{lem:phew!} implies that
\[\text{distance}(\lambda_v,\sigma(S))\le\|(S-\lambda_v)\varphi_v\|_2.\]

But, from Theorem~\ref{thm:wow!}, we also have for $i\in L$ that $(S\varphi_v-\lambda_v\varphi_v)(i)=\varphi_v(i)\cdot\big(\ell^*(v,i)-\lambda_v\big)$. As a result
\begin{align}
\nonumber\|(S-\lambda_v)\varphi_v\|_2^2 
&=\sum_{i\in L(v)}\varphi_v^2(i)\cdot\big(\ell^*(v,i)-\lambda_v\big)^2\\
\label{ide:superuseful} &=\frac{\rho_{v1}}{|L(v0)|}\sum_{i\in L(v0)}(\ell^*(i,v)-\lambda_v)^2+\frac{\rho_{v0}}{|L(v1)|}\sum_{i\in L(v1)}(\ell^*(i,v)-\lambda_v)^2,
\end{align}
where for the last identity we have used the equation~(\ref{lem:avgOFavg}). 
To complete the proof of the theorem note that $(\rho_{v0}+\rho_{v1})=1$. In particular, from the identity in equation~(\ref{ide:lambda*}), we may rewrite
\begin{align*}
\sum_{i\in L(v0)}(\ell^*(i,v)-\lambda_v)^2
&=
\sum_{i\in L(v0)}\left(\rho_{v0}\Big\{\overline{\ell^*(L(v1),v)}-\overline{\ell^*(L(v0),v)}\Big\}+\ell^*(i,v)-\overline{\ell^*\big(L(v0),v\big)}\right)^2\\
&=|L(v0)|\,\rho_{v0}^2\,\Big\{\overline{\ell^*(L(v1),v)}-\overline{\ell^*(L(v0),v)}\Big\}^2\\
&\qquad\qquad\qquad\qquad\qquad\qquad\qquad\quad+\sum_{i\in L(v0)}\left(\ell^*(i,v)-\overline{\ell^*\big(L(v0),v\big)}\right)^2.
\end{align*}
Namely
\[\frac{1}{|L(v0)|}\sum_{i\in L(v0)}(\ell^*(i,v)-\lambda_v)^2=\rho_{v0}^2\Big\{\overline{\ell^*(L(v1),v)}-\overline{\ell^*(L(v0),v)}\Big\}^2+\var\big(\ell^*(L,v);L(v0)\big).\]
Similarly,
\[\frac{1}{|L(v1)|}\sum_{i\in L(v1)}(\ell^*(i,v)-\lambda_v)^2=\rho_{v1}^2\Big\{\overline{\ell^*(L(v1),v)}-\overline{\ell^*(L(v0),v)}\Big\}^2+\var\big(\ell^*(L,v);L(v1)\big).\]
The second identity in the theorem is now a direct consequence of~(\ref{ide:superuseful}) and the last two identities.

Finally, again due to Theorem~\ref{thm:wow!}, we find that
\[\cos(S\varphi_v,\lambda_v\varphi_v)
=\frac{\varphi_v'S\varphi_v}{\|S\varphi_v\|_2}
=\frac{\lambda_v}{\sqrt{\varphi_v'\diag(\ell^*(L,v))^2\varphi_v}}.\]
But, similarly as we argued before
\[\varphi_v'\diag(\ell^*(L,v))^2\varphi_v
=\sum\limits_{i\in L}\varphi_v^2(i)\ell^*(i,v)^2
=\lambda_v^2+\sum\limits_{i\in L}\varphi_v^2(i)(\ell^*(i,v)-\lambda_v)^2
=\lambda_v^2+\|(S-\lambda_v)\varphi_v\|_2^2,\]
which implies that
\[\cos(S\varphi_v,\lambda_v\varphi_v)=\frac{1}{\sqrt{1+\frac{\sum\limits_{i\in L}\varphi_v^2(i)(\ell^*(i,v)-\lambda_v)^2}{\lambda_v^2}}}.\]
The third identity in the theorem follows now from equation (\ref{ide:superuseful}).
\end{proof}

It follows from the theorem that $S\varphi=\lambda_v\varphi_v$ if and only if $\overline{\ell^*(L(v1),v)}=\overline{\ell^*(L(v0),v)}$ and $\var\big(\ell^*(L,v);L(v0)\big)=\var\big(\ell^*(L,v);L(v1)\big)=0$. But these conditions are precisely equivalent to having the ORB-tree trace-balanced at $v$. Lemma~\ref{lem:varphiveig} may be therefore regarded a corollary of~Theorem~\ref{thm:eigest}.

We also emphasize that the first upper-bound for $\text{distance}\big(\lambda_v,\sigma(S)\big)$ may be computed efficiently using Theorem~\ref{thm:wow!}. Nevertheless, its alternative expression gives a way to quantify how approximately trace-balanced an ORB-tree is.

\begin{figure}
\centering
\includegraphics[scale=0.5]{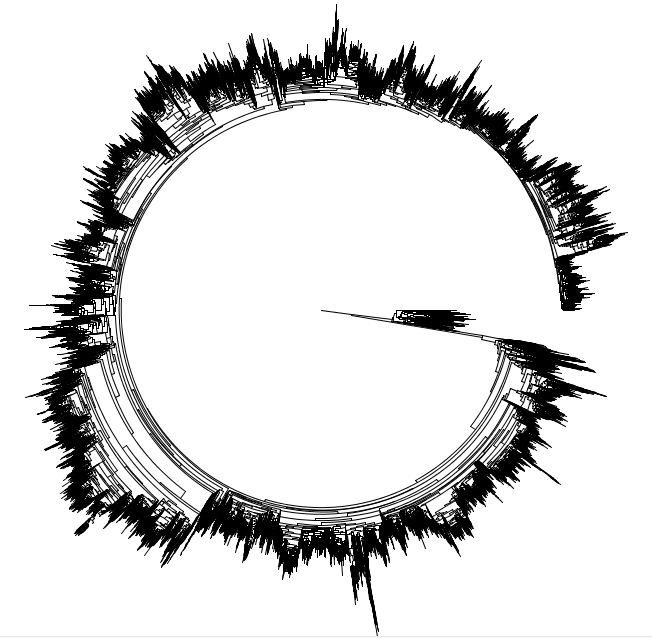}
\caption{\textbf{Circular layout and basic statistics of the 97\% Greengenes tree.} The tree has 99,322 leaves, 198,642 edges, and height (i.e. maximal leaf depth) 107. The average branch length is $1.42\times10^{-2}$ units, with lengths varying between $1.5\times10^{-4}$ and $1.0$. }
\label{fig:circular97greengenes}
\end{figure}

\section{New Insights into the Microbial Tree of Life}
\label{sec:treeoflife} 

Many methods in microbiology rely on a phylogenetic tree relating microorganisms. At the microbial level, however, the notions of genus or species are ill-defined because microorganisms do not interbreed. So microbes' taxonomy and phylogeny are often based on the so-called 16S ribosomal RNA (16S rRNA) gene. This gene is present in all known single cell organisms and can therefore be used as a phylogenetic marker. An operational taxonomic unit (OTU) is a cluster of these markers defined by some least level of DNA sequence similarity among its (highly) conserved regions.

Greengenes is a standarized database based on the 16S rRNA marker. It has been a standard reference in microbial studies, particularly metagenomics, and is the default option in QIITA \cite{Qiita18}---a widely used open-source management platform for microbial analyses. Greengenes phylogenetic trees are built using FastTree \cite{prideh10} and their associated taxonomies are assigned using tax2tree \cite{mcdpri11}. Trees are typically stored in the newick format~\cite{FelArcDayetal86}, which encodes their topology and branch lengths. Trees can be visualized from their newick format using software such as FigTree \cite{FIG} or ETE Toolkit \cite{hue16}.

Figure~\ref{fig:circular97greengenes} displays the Greengenes tree when OTUs are thresholded at a 97\% sequence similarity---the average similarity of macro-organisms' DNA in the same species. The tree represents the inferred evolutionary history of modern day microorganisms from common ancestors. Its root is at the center of the circular layout, and each OTU is associated with a single leaf in the tree and vice versa. Branch lengths are a proxy of evolutionary time such as the estimated expected number of mutations per nucleotide site~\cite{hoduc14}, and interior nodes (also called splits) are inferred speciation events that have led to the present-day microorganisms in the database.

A fundamental problem in microbiology is to link environmental factors (such as acidity, light, nutrients, salinity, temperature, etc) with microbial composition. A valuable tool for this has been the concept of $\beta$-diversity (i.e., a measure of differences between microbial composition across different environments). Early approaches to $\beta$-diversity include the Bray-Curtis dissimilarity~\cite{BraCur57} or the Jaccard distance~\cite{JAC}, which ignored the evolutionary relationships between microorganisms found in different environments. Nonetheless, one would expect microbes with a shared evolutionary history to similarly thrive or struggle in similar environments. Phylogenetic informed metrics were introduced precisely to convey this idea. These metrics require a phylogenetic tree relating the microorganisms observed in samples from all the environments under study. Among other more recent phylogentic trees such as SILVA \cite{quast12} and WoL \cite{zhu19}, Greengenes has been a common choice of representative phylogeny. We base our application on the latter---though our discussion applies to any phylogeny.

Double Principal Coordinate Analysis (DPCoA)~\cite{pav04} is a phylogenetically informed $\beta$-diversity metric between pairs of microbial environments, which provides similar insights~\cite{FukMcMDetEtAl12} to other more recent though more widely used  distances such as unweighted and weighted UniFrac~\cite{LozKni05}.

Let $T$ be the ORB-tree associated with a phylogenetic tree (e.g. the 97\% Greengenes tree), and $S$ the covariance matrix of $T$. In the context of phylogenetic informed metrics, environments are represented as probability mass functions over the OTUs (i.e. leaves). We denote those functions with lower-case letters such as $a$ and $b$, and interpret them as probability models over $L$. In particular, $a:L\to[0,+\infty)$ satisfies that $\sum_{x\in L}a(x)=1$ and, for each $e\in E$, $a(e)=\sum_{x\in e}a(x)$. With this convention, the DPCoA distance between two environments $a$ and $b$ is defined as~\cite{pav04,FukMcMDetEtAl12}: 
\begin{equation}
d(a,b):=\left\{\sum_{e\in E}\ell(e)\,\big(a(e)-b(e)\big)^2\right\}^{1/2}=\sqrt{(a-b)'S(a-b)}.
\label{ide:Dab}
\end{equation}
In particular, since $S$ is positive definite, DPCoA corresponds to a Mahalanobis distance~\cite{mah36}; implying that $d(\cdot,\cdot)$ is a metric---in the mathematical sense---in $\mathbb{R}^{|L|}$.

The weighted and unweighted UniFrac distances are instead defined as follows~\cite{LozKni05}:
\begin{align*}
d_w(a,b) &:=\sum_{e\in E}\ell(e)\,|a(e)-b(e)|;\\
d_u(a,b) &:=\dfrac{\sum_{e\in E}\ell(e)\,\big|\indicator[a(e)>0]-\indicator[b(e)>0]\big|}{\sum_{e\in E} \ell(e)}.
\end{align*}
Both versions of UniFrac are known to satisfy the triangular inequality~\cite{LozLlaKnietal11}. DPCoA is also more robust to unbiased noise but more sensitive to outliers than UniFrac~\cite{FukMcMDetEtAl12}.

Regardless of the metric of choice, the standard approach to linking environmental factors with microbial composition goes roughly as follows~\cite{LlaKni13}. First, environmental samples are collected, and each environment is represented by its OTU composition on the leaves of the phylogeny of reference. Then, the pairwise distance matrix between the environments is computed, and the environments are embedded into a low-dimensional Euclidean space using standard techniques such as multidimensional scaling (MDS)~\cite{BorGro05}. Despite the noisy and high-dimensional nature of microbial datasets~\cite{SogMorHubetal06,LlaGouRee11,HamLla12}, this approach has been remarkably reliable for the ordination~\cite{jon95} of microbial environments in as little as 1-2 dimensions, and for correlating environmental factors with microorganisms. However, this approach does not usually explain correlations, which need to be justified by other means.

In what remains of this section, we apply our methods to the Greengenes phylogeny. First, we demonstrate significant sparsification of the associated covariance matrix after applying the Haar-like wavelet transform. Then, we motivate a new wavelet-based phylogenetic $\beta$-diversity metric corresponding to a multiscale analysis of the phylogenetic tree. Finally, applying the new metric to a previously studied dataset shows that the wavelet-based metric can give novel insights into the relationship between environmental factors and OTU composition.

\subsection{Greengenes Phylogenetic Covariance Matrix Sparsification} The 97\% Greengenes tree has about 100,000 leaves. We can think of it as an ORB-tree by adding an external root $\circ$ and connecting it to the original root with a branch of length $0$. (Alternatively, we could think of the Greengenes tree as two ORB-trees with their roots merged.) We denote the resulting ORB-tree as $T$.

\begin{figure}
    \centering
    \begin{subfigure}{.49\textwidth}
    \centering
    \includegraphics[scale=.35]{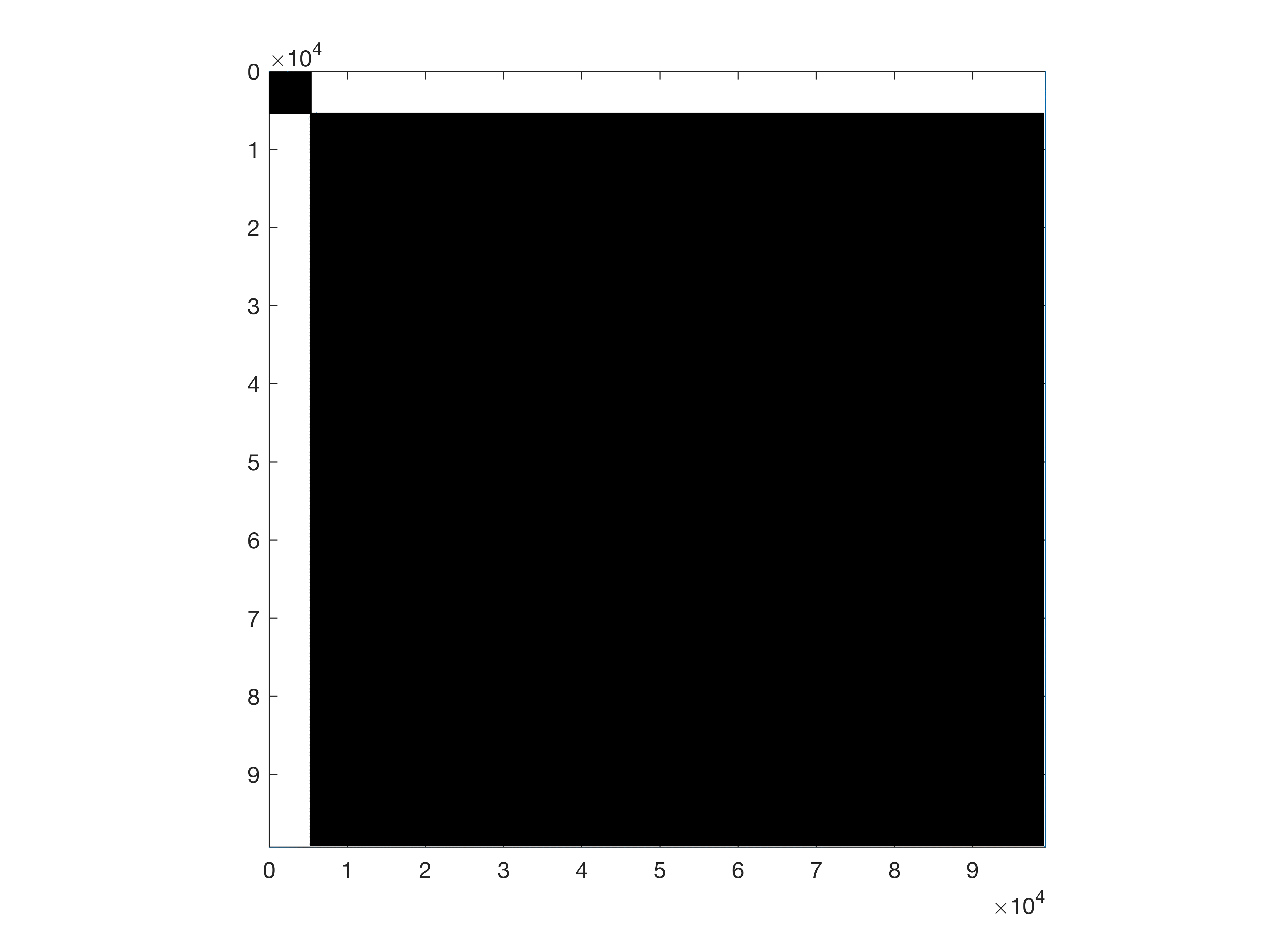}
    \caption{}
    \end{subfigure}
    \begin{subfigure}{.49\textwidth}
    \centering
    \includegraphics[scale=.25]{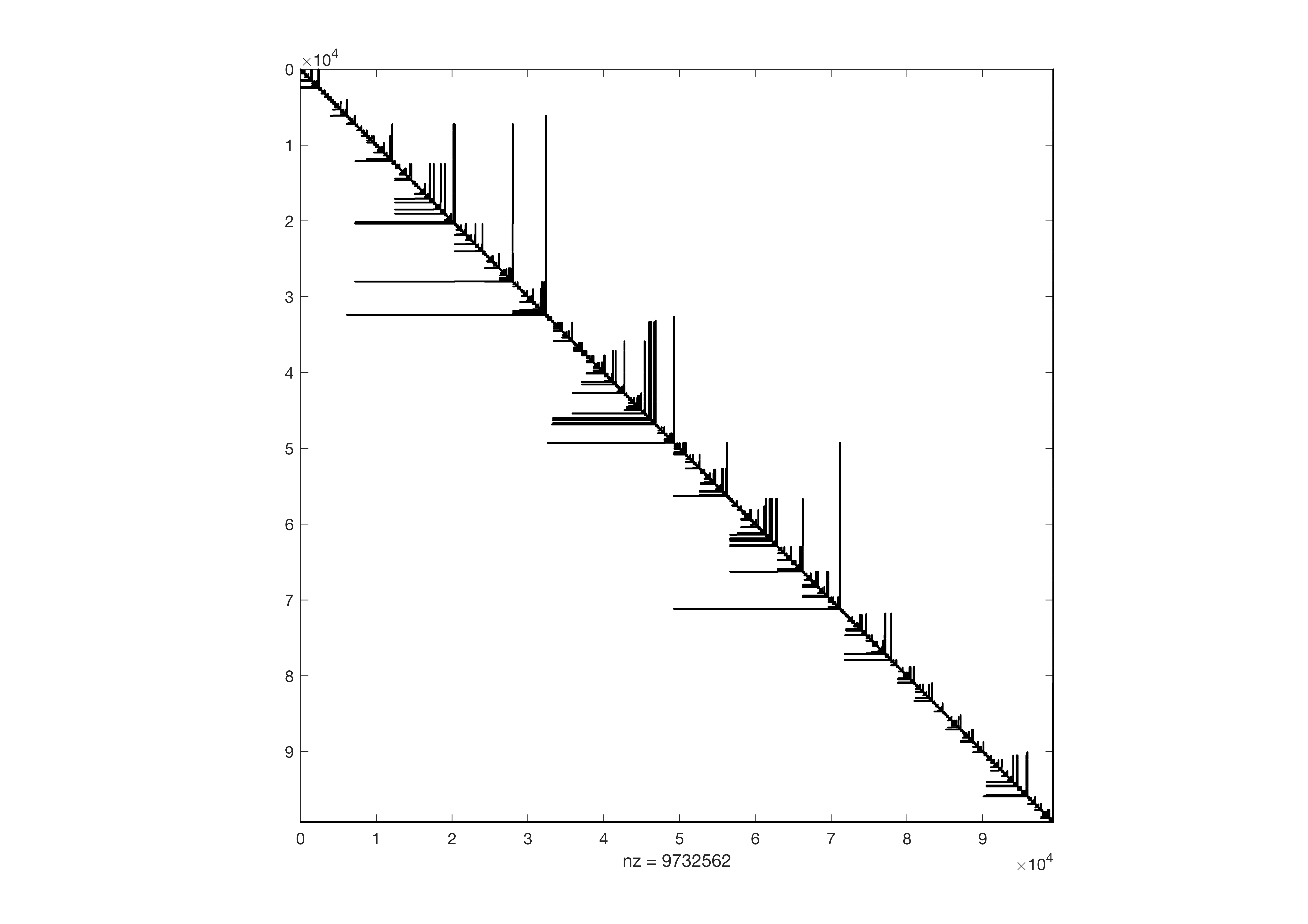}
    \caption{}
    \end{subfigure}
    \caption{\textbf{Heatmaps of matrices associated with the 97\% Greengenes.} Black (white) pixels denote non-zero (vanishing) entries. (a) Phylogenetic covariance matrix $S$ of the 97\% Greengenes tree. $S$ has dimensions $\sim 10^5\times10^5$. (b) Sparsified matrix $\Phi'S\Phi$.}
    \label{fig:97sparsification}
\end{figure}

The identity in equation~(\ref{ide:Sij}) implies that the covariance matrix $S$ of $T$ is a $2\times 2$ block diagonal matrix, with each block corresponding to an ORB-subtree. Nevertheless, approximately 94\% of the almost 10 billion entries in $S$ are non-zero because one of the ORB-subtrees (corresponding to the Archaea domain) is much smaller than the other---see Figure~\ref{fig:97sparsification}(a). This makes storing the covariance matrix of $T$ challenging. Further, basic computational tasks such as finding the spectrum and inverting $S$ for parameter estimation in phylogenetic comparative methods \cite{jhw22,fre12} is infeasible because this large matrix is almost fully dense. We may use, however, the Haar-like matrix $\Phi$ associated with $T$ to sparsify $S$. From Theorem~\ref{thm:zetaexplicit}, we can guarantee that $\zeta\ge0.9989$, i.e. at least 99.89\% of the entries in the similar matrix $\Phi'S\Phi$ vanish. This significant compression of the matrix $S$ can be appreciated in Figure~\ref{fig:97sparsification}(b).

\begin{figure}
    \centering
    \includegraphics[scale=.75]{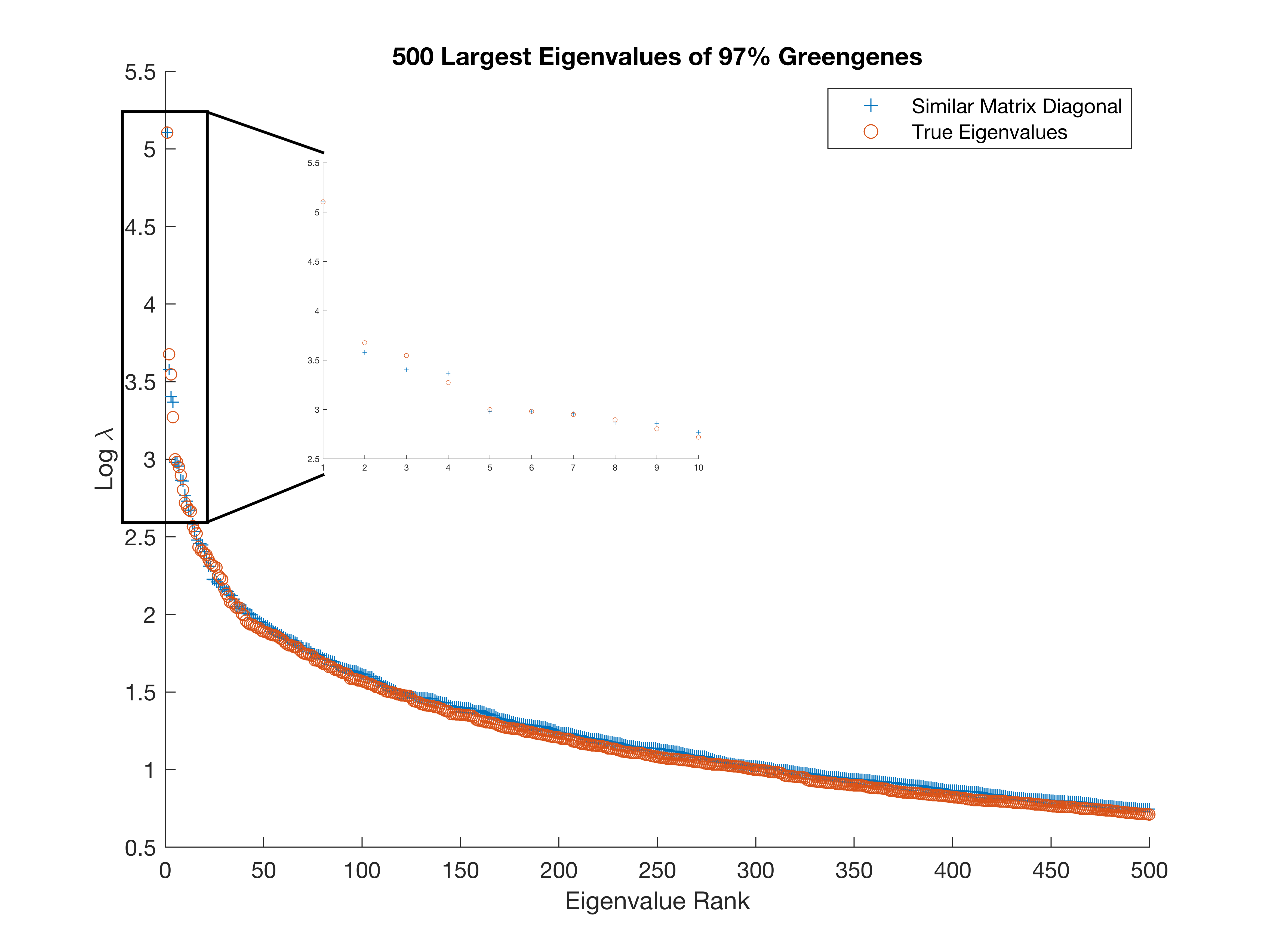}
    \caption{\textbf{Spectrum decay of 97\% Greengenes tree covariance matrix and corresponding approximation using Haar-like wavelets.} Only the 500 most dominant eigenvalues of $S$ are plotted as a function of their rank. Logarithms are in base-10.}
    \label{fig:circular97approxBYrank}
\end{figure}

We implemented Algorithm~\ref{alg:sparsify} using the sparse matrix packages from SciPy \cite{2020SciPy-NMeth} to compute $\Phi'S\Phi$. As proof-of-principle we used this compressed representation of $S$ to compute its largest 500 eigenvalues to machine precision using SciPy's implementation of the Lanczos algorithm. As seen in Figure~\ref{fig:circular97approxBYrank}, the eigenvalues of $S$ decay rapidly. In fact, we found that $\lambda_1(S)\sim 1.27\times10^5$, $\lambda_2(S)\sim 4.75\times10^3$, and $\text{trace}(S)\sim 1.65\times 10^5$, so the top and top-two eigenvalues account for approximately 77\% and 80\% of the trace of $S$, respectively.

As seen in Figure~\ref{fig:circular97approxBYrank} also, the sorted diagonal entries in $\Phi'S\Phi$ (i.e. the quantities $\lambda_v$, with $v\in I$) approximate with ample accuracy the spectrum of $S$. For instance, $\max_{v\in I}\lambda_v$ underestimates $\lambda_1(S)$ with only about a $0.06\%$ relative error. Anticipating this overall accuracy from $T$ alone remains an open problem as neither our mathematical results, particularly Theorem~\ref{thm:eigest}, nor more general ones such as the Gershgorin's circle theorem, Sylvester's determinant theorem, and bounds found in \cite{Tho76,bunnie78,ipsnad09} have been able to explain it.

\subsection{A Wavelet Based Phylogenetic $\beta$-diversity Metric}
\label{sec:Haardist}

Let $T$ be the ORB-tree associated with a phylogenetic tree. Recall that $\varphi_v$, with $v\in I$, is supported on $L(v)$, and together these functions form an orthonormal base of $\mathbb{R}^{|L|}$. In particular, just as wavelets are traditionally used to localize signals at different scales, we may use the Haar-like basis of $T$ to localize environmental OTU distributions on subsets of leaves defined by splits in the tree. This is particularly appealing from a biological standpoint. Indeed, the opposite signs of $\varphi_v$ on the leaves of the left and right subtrees dangling from $v$ may be interpreted as a speciation event that conferred more fitness to present-day microorganisms descending from one of the subtrees than the other. We propose the following definition to convey these features into a phylogenetic $\beta$-diversity metric.

Recall that $\lambda_v=(\Phi'S\Phi)(v,v)>0$, for each $v\in I$. Further, for a given environment $a$ (i.e., OTU distribution over $L$), $\Phi'a$ is the projection of $a$ onto the Haar-like basis of the reference tree.

\begin{tcolorbox}
\begin{definition}
The Haar-like distance between two environments $a$ and $b$ is the quantity
\[d_h(a,b):=\sqrt{\sum_{v\in I}\lambda_v\,\Delta_v^2},\text{ where }\Delta=(\Delta_v)_{v\in I}:=\Phi'(a-b).\]
\end{definition}
\end{tcolorbox}

The specifics of this distance can be motivated as follows. On one hand, the terms $\Delta_v^2$, with $v\in I$, convey the idea that $d_h$ regards two environments similar (different) when their OTU compositions project similarly (differently) onto the Haar-like basis of the reference tree. On the other hand, the weights $\lambda_v$, with $v\in I$, are motivated by the success of DPCoA in various biological investigations. To explain this, consider the matrices $D:=\text{diag}(\lambda_v:v\in I)$ and $E:=\Phi'S\Phi-D$. Observe that $d_h(a,b)=\sqrt{\Delta'D\Delta}$; in particular, $d_h$ is a metric in $\mathbb{R}^{|L|}$ because $D$ is positive definite, and $d(a,b)=\sqrt{\Delta'D\Delta+\Delta'E\Delta}$. In large phylogenetic trees, however, we expect $E$ to be mostly filled with zeroes due to Corollary~\ref{cor:randomT}---which suggests considering $d_h$ as an alternative metric to DPCoA.

We have mentioned before that while traditional phylogenetic metrics (in conjunction with embedding techniques) have been remarkably successful at correlating microbial composition with environmental factors, these correlations cannot usually be explained from the metrics alone. The wavelet nature of the Haar-like distance has, however, the potential to explain said correlations. Indeed, the biological interpretation of the Haar-like basis conveyed by their sign flip suggests that if $\lambda_v\,\Delta_v^2$ is comparatively large (small) for some $v\in I$, then the speciation event associated with $v$ has a significant (little) influence differentiating the OTU distributions between two environments $a$ and $b$. (There may be discrepancies between taxonomy and splits in a tree. In particular, while the aforementioned correlations may be explained by a phylogeny, they are not necessarily explained by a taxonomic classification.)

\subsection{Haar-like Distances of the Guerrero Negro microbial mat} 

A microbial mat is a bio-film of layered groups of microorganisms with coupled biochemistries. Their rich biodiversity, combined with the environmental gradients of light, oxygen, etc., offer an ideal setting to test phylogenetic $\beta$-diversity metrics.

To demonstrate the insight gained from Haar-like distance we applied it to a 16S rRNA data set of 18 soil samples (obtained from QIITA~\cite{Qiita18}) at different depths of the Guerrero Negro microbial mat~\cite{kir12}, located in Baja California Sur, Mexico. We used the 97\% Greengenes as the reference phylogeny.

\begin{figure}
\centering
\includegraphics[scale=.4]{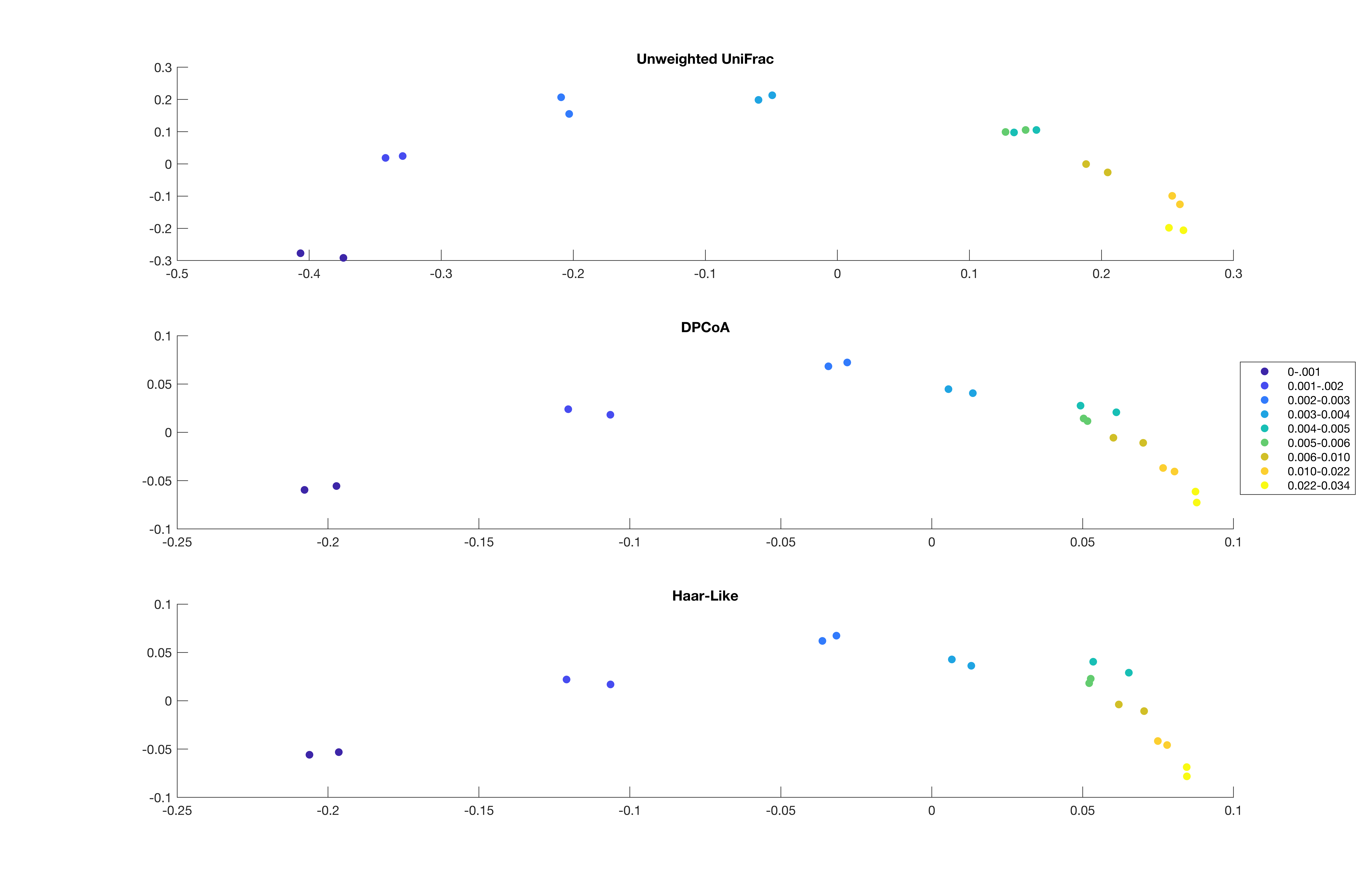}
\caption{\textbf{2-D MDS embeddings of samples from Guerrero Negro w.r.t. different metrics.} The embeddings are based on unweighted UniFrac (top), DPCoA (middle), and the Haar-like distance (bottom). Depth varies from 0-0.034 meters.}
\label{fig:microbialmat}
\end{figure}

Earlier work~\cite{kir12} based on unweighted UniFrac showed a gradient of microbial composition in the mat with respect to depth---see top plot in Figure~\ref{fig:microbialmat}. (For a discussion regarding the ``horseshoe" shape in the plot see \cite{kir12,mor17,dia08}.) As seen on the bottom plot of the same figure, we can practically reproduce this gradient using the Haar-like distance instead. Furthermore, as seen on the bottom two plots, the DPCoA and Haar-like distance produce nearly indistinguishable embeddings.

\begin{figure}
\centering
\includegraphics[scale=.6]{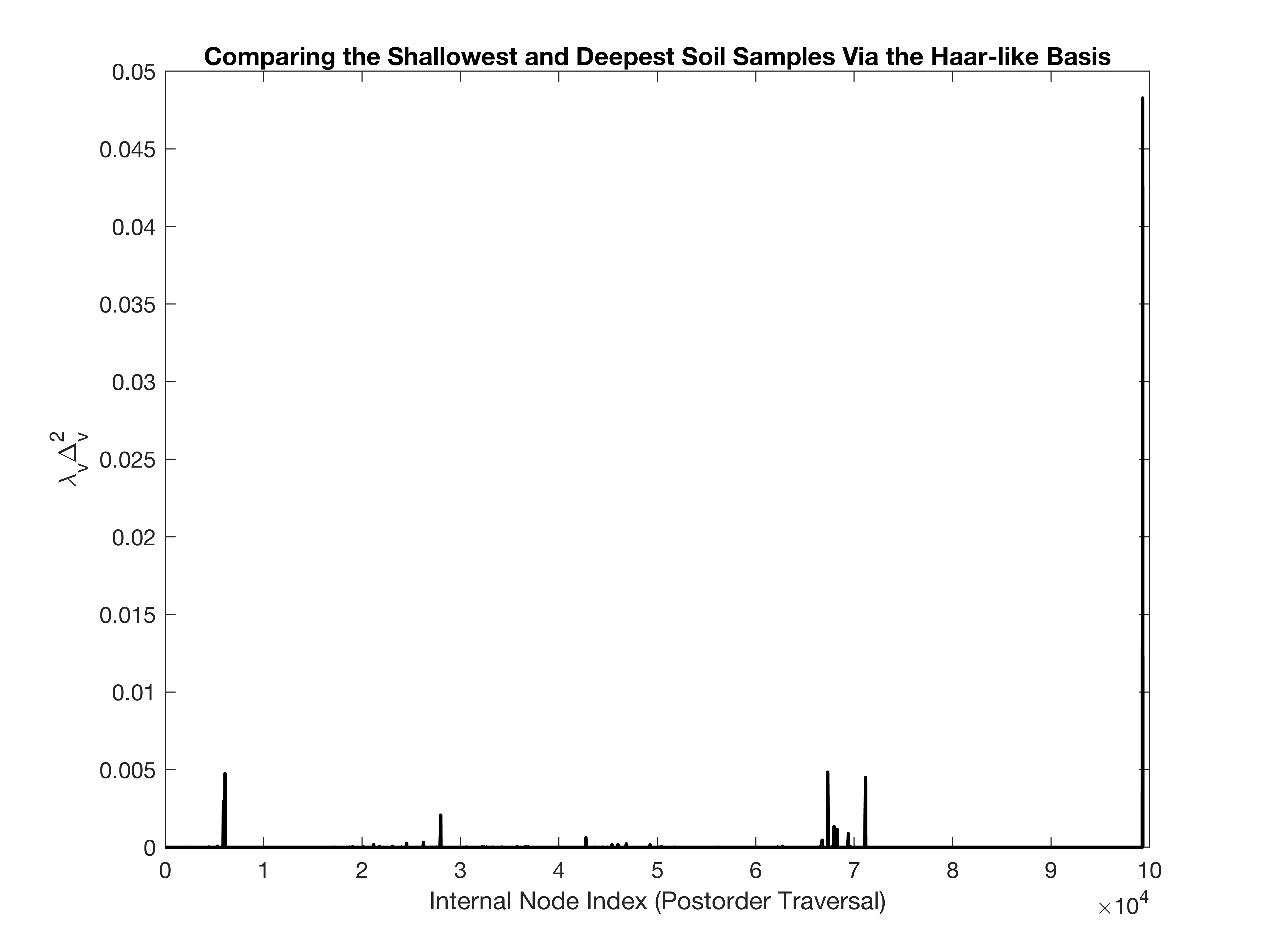}
\caption{\textbf{Plot of $\mathbf{\lambda_v\,\Delta_v^2}$, with $v\in I$, to measure the Haar-like distance between the shallowest and deepest sample in the Guerrero Negro dataset.} The average non-zero value of $\lambda_v\,\Delta_v^2$ is $\sim 2.09\times 10^{-5}$. The standard deviation of these values is $\sim 7.55\times 10^{-4}$.}
\label{fig:microbialmatproj}
\end{figure}

\begin{figure}
\centering
\includegraphics[width=\linewidth]{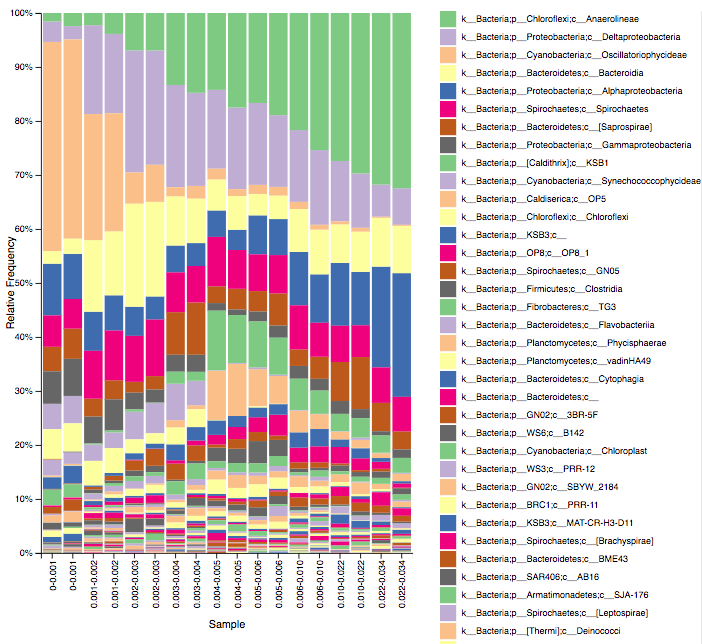}
\caption{\textbf{OTU compositions plotted at class level taxonomy in QIITA.} Only the 35 most abundant OTU class names are displayed in the legend.}
\label{fig:microbialmatcomp}
\end{figure}

While the three phylogenetic $\beta$-diversity metrics imply that soil depth drives a measurable change in OTU composition, we can go a step further with the Haar-like distance and determine which splits are responsible for this trend and quantify their importance. We demonstrate this by comparing the two extremes in the dataset: let $a$ and $b$ be the OTU compositions of the shallowest and deepest environment, respectively. Define $\Delta=\Phi'(b-a)$. Following the logic described in Section~\ref{sec:Haardist}, we computed $v\in I\longrightarrow\lambda_v\,\Delta_v^2$, indexing interior nodes according to a postorder traversal of the 97\% Greengenes tree. As seen in Figure~\ref{fig:microbialmatproj}, the 3 largest values are statistically significant. These are associated with the Haar-like wavelets $\varphi_{99311}$ ($\lambda_v\Delta_v^2$-value $\sim 4.84\times10^{-2}$), $\varphi_{67317}$ ($\lambda_v\Delta_v^2$-value $\sim 4.84\times10^{-3}$), and $\varphi_{6079}$ ($\lambda_v\Delta_v^2$-value $\sim 4.75\times10^{-3}$). These correspond to splits at depths 10, 34, and 18 of the 97\% Greengenes tree, respectively.

Notably, the split associated with $\varphi_{99311}$ corresponds to the largest $\lambda_v\Delta_v^2$-value. According to the associated taxonomic classification, the (say) left descendants of this split associated correspond to the phylum level classification of Cyanobacteria. This is consistent with the conclusion in~\cite{kir12}, which correlated Cyanobacteria abundance changes with soil depth and explained the correlation by their ability to photosynthesize.

The other two wavelets provide novel insight into other important OTU composition differences driving the observed soil depth gradient in the Guerrero Negro dataset. Indeed, while the descendants of the split associated with $\varphi_{67317}$ do not exhaust a taxonomic classification, all leaves under that split are classified as Anaerolineae. This differentiation between the shallowest and deepest sample may be due to Anaerolineae's role as a anaerobic digester~\cite{xia16}. This claim is reinforced by the top green bars in the OTU composition plot in Figure~\ref{fig:microbialmatcomp}, which show a rapid increase of Anaerolineae with depth. 

Finally, the split associated with $\varphi_{6079}$ subdivides the Cyanobacteria phylum into further classes, including Oscillatoriophycideae which, according to Figure~\ref{fig:microbialmatcomp}, is the third most abundant class of the Guerrero Negro dataset. The relevance of this split to differentiate shallow from deep samples may be explained by Oscillatoriophycideae's photoautotrophic capability~\cite{StaCoh77}. Again, this is supported by Figure~\ref{fig:microbialmatcomp}, which shows a sharp decrease in Oscillatoriophycideae with respect to soil depth.

Our analysis of the Guerrero Negro mat shows that the Haar-like distance may be a valid alternative to other more common phylogenetic $\beta$-diversity metrics, primarily because it provides a systematic method for detecting statistically significant speciation events (and corresponding levels of OTU classification) that can link OTU composition with environmental factor gradients.

\appendix
\section{Orthonormality of Haar-like bases} 

The statement that the Haar-like basis $\{\varphi_v\}_{v\in I}$ associated with an ORB-tree is orthonormal is based on the concept of multiresolution analysis of Euclidean spaces in~\cite{gavnad10}. Here we justify this fact by first principles.

Let $u,v\in I$. If $u=v$ then
\[\langle \varphi_u,\varphi_v\rangle=\frac{|L(u1)|+|L(u0)|}{|L(u)|}=1.\]
Instead, there are two possibilities when $u\ne v$. If $L(u)\cap L(v)=\emptyset$ then $\langle \varphi_v,\varphi_u\rangle=0$ because $\varphi_v$ and $\varphi_u$ have disjoint supports. Otherwise, if $L(u)\cap L(v)\ne\emptyset$ then Lemma~\ref{lem:LuLvnotempty} let us assume without any loss of generality that $u$ is an ancestor of $v$. In particular, $L(v)\subset L(u)$ but also $\varphi_u$ remains constant over $L(v)$. Therefore, for any given $x\in L(v)$:
\begin{align*}
\langle \varphi_u,\varphi_v\rangle
&=\varphi_u(x)\cdot\sum_{y\in L(v)}\varphi_v(y)\\
&=\varphi_u(x)\cdot\left\{\sqrt{\frac{ |L(v1)|\cdot|L(v0)|}{|L(v)|}}-\sqrt{\frac{ L(v0)|\cdot|L(v1)|}{|L(v)|}}\right\}=0.
\end{align*}

\color{black}

\section*{Acknowledgments}
This work has been partially funded by the NSF grant No. 1836914.

\bibliographystyle{siamplain}
\bibliography{references}

\end{document}